\documentclass[12pt]{report} 
\usepackage{latexsym}
\newcommand{\hs}{\hspace{0.3cm}}

\newcommand{\rank}{\mbox{rank}}

\newcommand{\diag}{\mbox{diag}}

\newcommand{\IZ}{Z\!\!\!Z}

\newcommand{\IC}{I\!\!\!\!C}
\newcommand{\IQ}{I\!\!\!\!Q}

\newcommand{{\CB}}{\cal B}

\newcommand{{\CC}}{\cal C}

\font\myfonts=cmfi10 scaled \magstep 1
 1

\def\picture#1 by #2 (#3){
\vbox to #2{
\hrule width #1 height 0pt depth 0pt
\vfill
\special{picture #3}}}

\begin{document}
\begin{center}
{\bf \Large Computing invariants and semi-invariants by means of Frobenius Lie algebras}\\
\ \\
{\bf \large Alfons I. Ooms}\\
\ \\
{\it Mathematics Department, Hasselt University,\\
Agoralaan, 3590 Diepenbeek, Belgium\\
E-mail address: alfons.ooms@uhasselt.be}
\end{center}
\ \\
\ \\
\pagestyle{plain}
{\bf Abstract.}\hs Let $U(${\myfonts g}$)$ be the enveloping algebra of a finite dimensional Lie algebra {\myfonts g} over a field $k$ of characteristic zero, $Z(U(${\myfonts g}$))$ its center and $Sz(U(${\myfonts g}$))$ its semicenter.  A sufficient condition is given in order for $Sz(U(${\myfonts g}$))$ to be a polynomial algebra over $k$. Surprisingly, this condition holds for many Lie algebras, especially among those for which the radical is nilpotent, in which case $Sz(U(${\myfonts g}$)) = Z(U(${\myfonts g}$))$.  In particular, it allows the explicit description of $Z(U(${\myfonts g}$))$ for more than half of all complex, indecomposable nilpotent Lie algebras of dimension at most 7.\\
\ \\
\begin{center}
{\bf 1. INTRODUCTION}
\end{center}
Let $U(${\myfonts g}$)$ be the enveloping algebra of a nonzero finite dimensional Lie algebra {\myfonts g} over a field $k$ of characteristic zero, $Z(U(${\myfonts g}$))$ its center and $Sz(U(${\myfonts g}$))$ its semicenter, i.e. the subalgebra of $U(${\myfonts g}$)$ generated by the semi-invariants of  $U(${\myfonts g}$)$.  The (Poisson) semicenter $Sz(S(${\myfonts g}$))$ of the symmetric algebra $S(${\myfonts g}$)$, which is defined in a similar fashion, is known to be isomorphic to $Sz(U(${\myfonts g}$))$ if $k$ is algebraically closed $[RV]$.  Unlike the center $Z(U(${\myfonts g}$))$ [DNOW, p. 321], the semicenter $Sz(U(${\myfonts g}$))$ is never reduced to $k$ [D1, 4.9.4] and is always factorial [Moe], [LO].\\
In some important cases $Sz(U(${\myfonts g}$))$ turns out to be a polynomial algebra over $k$.  This happens for instance if {\myfonts g} is: semisimple [D1, 7.3.8], Frobenius (i.e. the index $i(${\myfonts g}$) = 0$) [DNO, p. 339], square integrable (i.e. $i(${\myfonts g}$) = dim\ Z(L)$) [DNOW, p. 323], the Lie algebra of strictly lower triangular matrices [D4], any Borel subalgebra (as well as its nilradical) of a complex semisimple Lie algebra [J2, p. 263 and p. 257] and many (but unfortunately not all [Y]) parabolic (and biparabolic) subalgebras of certain complex semisimple Lie algebras [FJ1,2], [J3,4].\\
\ \\
The key result of this paper is the following.  Let {\myfonts g} be a Lie algebra for which there exists a split torus \quad $T \subset Der\ ${\myfonts g}\quad (i.e. a commutative Lie algebra consisting of diagonalizable derivations of {\myfonts g}) such that the semi-direct product \quad $L = T\ \oplus\ ${\myfonts g}\quad is a Frobenius Lie algebra. Let \quad $x_1,\ldots,x_n$\quad be a basis of $L$ and let \quad $v_1,\ldots, v_r$\quad be the irreducible factors of the nonzero determinant \quad $det([x_i,x_j]) \in S(${\myfonts g}$)$. \quad Then \quad $v_1,\ldots,v_r$ \quad are the only (up to nonzero scalar multipliers) irreducible semi-invariants of $S(L)$ and their weights \quad $\lambda_1,\ldots, \lambda_r \in L^\ast$ \quad are linearly independent over $k$ [DNO, pp. 339-343]. We will show that $S(${\myfonts g}$)$ and $S(L)$ have the same semicenter, which coincides with the polynomial algebra \quad $k[v_1,\ldots, v_r]$.\quad
A similar result holds for the semicenters of $U(${\myfonts g}$)$ and $U(L)$.  Under these circumstances \quad $dim\ T = i(${\myfonts g}$)$,\quad which implies that \quad $i(${\myfonts g}$) \leq rank$ {\myfonts g}, \quad the latter being the dimension of a maximal torus inside $Der$ {\myfonts g}.  Also, the sum of the degrees of \quad $v_1,\ldots, v_r$\quad is at most \\ 
$c(${\myfonts g}$)= (dim$ {\myfonts g}$\ +\ i(${\myfonts g}$))/2$,\quad an inequality which is known to be valid in the cases mentioned above and which has recently been shown to hold in general [OV].\\
In section 4 we exhibit examples of different types where this simple method can be employed successfully in order to compute algebraically independent generators for $Sz(U(${\myfonts g}$))$.\\
Section 5 is devoted to complex indecomposable nilpotent Lie algebras of dimension at most seven.  We use the same notation as Carles [Ca] and Magnin[Ma2], the latter providing, among other things, a maximal torus in each case.  However, we also include the notation of Morozov $(M_{\_})$ [Mor] and Cerezo $(C_{\_})$ [C2] in dimension 6 and of Seeley $((\ldots))$ [Se] and Romdhani $(R_{\_})$ [R] in dimension 7.  We distinguish two classes:\\
1) $i(${\myfonts g}$) \leq rank$ {\myfonts g}\\
Remarkably, in this situation our condition is satisfied with very few exceptions (namely, in dimension 6: {\myfonts g}$_{6.15}$ and in dimension 7: {\myfonts g}$_{7.3.1(i_\lambda),\lambda=0,1}$; {\myfonts g}$_{7,3.12}$; {\myfonts g}$_{7,3.13}$ and {\myfonts g}$_{7,3.17}$.  For all these {\myfonts g} is coregular, i.e. $Z(U(${\myfonts g}$))$ is a polynomial algebra, except for {\myfonts g}$_{7,3.17}$).\\
It turns out that the generators we find for $Z(U(${\myfonts g}$))$ coincide with the algebraic independent generators of its quotient field listed in [C2] for dimension 6 and in [R] for dimension 7 (although in the latter quite a few cases are missing [Go, pp. 146-149] and some minor corrections have to be made).  Most of the calculations were done with Maple.\\
As a welcome by-product we obtain all solvable, almost algebraic Frobenius Lie algebras over\ \ ${\IC}$ of which the nilradical is indecomposable of dimension at most seven. In each case a split torus $T \subset Der$ {\myfonts g} for which \quad $L = T \ \oplus\ ${\myfonts g} \quad is Frobenius is given. Sometimes, the possible solutions for $T$ form a family depending on parameters.  See also [CV].\\
2) $i(${\myfonts g}$) > rank$ {\myfonts g}\\
In this situation, the following methods are employed in order to find the generators of $S(${\myfonts g}$)^{\myfonts g}$ (and hence of $Z(U(${\myfonts g}$)))$: the procedure Dixmier used for dimension at most 5 [D3, pp. 322-330], Singular [GPS] (in case {\myfonts g} has an abelian ideal of codimension one), a result by Panyushev [P2, Theorem 1.2] and an effective algorithm recently constructed by Michel Van den Bergh.\\
We confine ourselves, for now, to finish the seven remaining cases in dimension 6, of which only two are not coregular, namely {\myfonts g}$_{6,6}$ and {\myfonts g}$_{6,16}$, the latter being the 6-dimensional standard filiform Lie algebra $L(6)$ (also called generic filiform).  In fact, the $n$-dimensional standard filiform Lie algebra $L(n)$ (i.e. with basis \quad $x_1,\ldots, x_n$ \quad and nonzero Lie brackets \quad $[x_1, x_i] = x_{i+1}$,\quad $i:2,\ldots, n-1$) is not coregular if and only if $n \geq 5$ [OV].

\begin{center}
{\bf 2. PRELIMINARIES AND GENERAL RESULTS}
\end{center}
{\bf (i) The index \boldmath{$i(L)$}}\\
Let $L$ be a Lie algebra over a field $k$ of characteristic zero with basis \quad $x_1,\ldots, x_n$.  Let \quad $f \in L^\ast$ \quad and consider the alternating bilinear form $B_f$ on $L$ sending \quad $(x,y)$ into \quad $f([x,y])$,\quad with kernel:
$$L(f) = \{x \in L\mid f([x,y]) = 0\ \mbox{for all}\ y \in L\}$$
We put \quad $i(L) = \min\limits_{f \in L^\ast}$ dim $L(f)$, the index of $L$.  We recall from [D1, 1.14.13] that
$$i(L) = \dim L - \rank_{R(L)}([x_i, x_j])$$
where $R(L)$ is the quotient field of the symmetric algebra $S(L)$ of $L$.  In particular $\dim L - i(L)$ is an even number.  Furthermore, $f$ is called regular if\quad $\dim L(f) = i(L)$.\quad It is well-known that the set $L^\ast_{\mbox{\tiny reg}}$ of all regular elements of $L^\ast$ is an open dense subset of $L^\ast$ for the Zariski topology.\\
Denote by $Z(D(L))$ the center of the quotient division ring $D(L)$ of $U(L)$ and let $R(L)^L$ be the subfield of the invariants of $R(L)$ under the action of $ad\ x$, $x \in L$.\\
\ \\
{\bf THEOREM 2.1.}
$$\mbox{trdeg}_k Z(D(L)) = \mbox{trdeg}_kR(L)^L \leq i(L)$$
[O3], [RV]. Moreover, equality occurs if one of the following conditions is satisfied:
\begin{itemize}
\item[(1)] $L$ is algebraic \quad [O3], [RV]
\item[(2)] $k$ is algebraically closed and $L$ has no proper semi-invariants in $S(L)$ (or equivalently in $U(L)$) \quad [OV, Proposition 4.1]
\end{itemize}
{\bf (ii) The semicenter \boldmath{$Sz(U(L))$}}\\
Let \quad $\lambda \in L^\ast$. We denote by $U(L)_{\lambda}$ the set of all \quad $u \in U(L)$ \quad such that \quad $ad\ x(u) = \lambda(x)u$\quad for all $x \in L$.  Any nonzero element \quad $u \in U(L)_\lambda$ \quad is called a semi-invariant w.r.t. the weight $\lambda$. We call $u$ a proper semi-invariant if $\lambda \neq 0$.  Clearly, \quad $U(L)_\lambda U(L)_\mu \subset U(L)_{\lambda + \mu}$\quad for all $\lambda, \mu \in L^\ast$.  The sum of the $U(L)_\lambda$ is direct and it is a subalgebra $
Sz(U(L))$ of $U(L)$, the semicenter of $U(L)$.\quad $Sz(U(L))$ is factorial and any semi-invariant can be written uniquely as a product of irreducible semi-invariants.  Suppose $z \in D(L), z \neq 0$.  Then, \quad $z \in Z(D(L))$\quad if and only if $z$ can be written as a quotient of two semi-invariants of $U(L)$ with the same weight.  So, if $L$ has no proper semi-invariants (as it is if its radical is nilpotent) then $Z(D(L))$ is the quotient field of $Z(U(L))$. A useful link with $S(L)$ is the symmetrization map, i.e. the canonical linear isomorphism $s$ of $S(L)$ onto $U(L)$, which maps each product \quad $y_1,\ldots,y_m,\ y_i \in L$, into \quad $(1/m!) \sum\limits_p y_{p(1)}\ldots y_{p(m)}$, where $p$ ranges over all permutations of $\{1,\ldots, m\}$.\\
$s$ is known to commute with each derivation of $L$ and hence maps $S(L)_\lambda$ onto $U(L)_\lambda$ and also $Sz(S(L))$ onto $Sz(U(L))$.  This restriction is an algebra isomorphism in case $L$ is nilpotent [D1, 4.8.12] or if $L$ is Frobenius [DNO, p. 399].\\
The weights of the semi-invariants of $U(L)$ (or $S(L)$) form a semigroup $\Lambda(L)$, which is not necessarily finitely generated [DDV].  However, the additive subgroup $\Lambda_D(L)$ of $L^\ast$ generated by $\Lambda(L)$, is a finitely generated free abelian group [NO].  Perhaps unaware of this a new proof was recently given [FJ, p. 1519].\\
The centralizer of $Sz(U(L))$ in $U(L)$ coincides with $U(L_\Lambda)$, where $L_\Lambda$ is the intersection of $\mbox{ker}\ \lambda$, $\lambda \in \Lambda(L)$.  $L_\Lambda$ is a characteristic ideal of $L$ for which the following hold [DNO, pp. 332, 343]:\\
\ \\
{\bf THEOREM 2.2.}
$$Z(U(L)) \subset Sz(U(L)) \subset Z(U(L_\Lambda)) = Sz(U(L_\Lambda))$$
and $Z(D(L_\Lambda))$ is the quotient field of $Z(U(L_\Lambda))$.  Moreover, \quad $Sz(U(L)) = Z(U(L_\Lambda))$\quad in case $k$ is algebraically closed and either $L$ is almost algebraic or $L$ is Frobenius.\\
Similar results hold in $S(L)$.  See also [FJ, p. 1518].\\
\ \\
{\bf (iii) Frobenius Lie algebras}\\
A Lie algebra $L$ is called Frobenius if there is a linear functional $f \in L^\ast$ such that the alternating bilinear form \quad $B_f(x,y) = f([x,y]),\ x,y \in L$, is nondegenerate, i.e. $i(L) = 0$.  This notion was introduced in [O1] in connection with the problem of Jacobson on the characterization of Lie algebras having a primitive universal enveloping algebra. See [O2,3] for the general solution. In particular, $U(L)$ is primitive if $L$ is Frobenius and the converse holds if $L$ is algebraic.  It has become clear over the years that these Lie algebras form a large class. They also play an important role in different areas.  For instance, many Lie subalgebras of semisimple Lie algebras, such as some parabolic and biparabolic (sea weed) subalgebras, are Frobenius [E1,2,3], [Ge], [P1], including most Borel subalgebras of simple Lie algebras [EO, p. 146].  Other examples are $j$-algebras which are real, solvable Lie algebras, which play an essential role in the study of bounded homogeneous domains in $\IC^n$, admitting a simply transitive solvable group of analytic automorphisms [PS].  Also, Frobenius (and quasi Frobenius) Lie algebras give rise to constant solutions of the classical Yang-Baxter equation [Dr].  More recently, Frobenius Lie algebras appear naturally in the study of simple hypersurface singularities [EK].\\
\\
We now collect some useful facts on semi-invariants from [O4], [DNO].  Let $L$ be a Frobenius Lie algebra with basis \quad $x_1,\ldots, x_n$.  Then $n$ is even and $L$ has a trivial center.  In fact, $L$ is highly noncommutative since \quad $Z(D(L)) = k$.  The Pfaffian \quad $Pf([x_i, x_j]) \in S(L)$\quad is homogeneous of degree $\frac{1}{2}\dim L$ and \quad $(Pf([x_i, x_j]))^2 = \det([x_i, x_j]) \neq 0$. \quad This is well determined by $L$ (up to a nonzero scalar multiple) and for all $g \in \mbox{Aut} L$:
$$g(Pf([x_i, x_j])) = \det g\ Pf([x_i, x_j])$$

{\bf THEOREM 2.3.}\quad Let $L$ be Frobenius and let \quad $v_1,\ldots, v_r \in S(L)$\quad be the (pairwise nonassociated) irreducible factors of $Pf([x_i, x_j])$.  Then,
\begin{itemize}
\item[(1)] $v_1,\ldots, v_r$ \quad are the only (up to nonzero scalar multipliers) irreducible semi-invariants of $S(L)$, say with weights $\quad \lambda_1,\ldots, \lambda_r \in \Lambda(L)$.
\item[(2)] $Sz(S(L)) = k[v_1,\ldots, v_r]$, \quad a polynomial algebra over $k$.
\item[(3)] $\lambda_1,\ldots, \lambda_r$ \quad are linearly independent over $k$.  They generate the semigroup $\Lambda(L)$ and \quad $L_\Lambda = \cap \ \mbox{ker}\  \lambda_i$, \quad $i : 1,\ldots, r$
\item[(4)] $r = \dim L - \dim L_\Lambda = i(L_\Lambda)$\quad (see also Lemma 3.1)
\end{itemize}

{\bf COROLLARY 2.4.} \quad Let $L$ be Frobenius.  Any semi-invariant of $S(L)$ is homogeneous.  It is also a semi-invariant under the action of $\mbox{Der}\ L$.\\
\ \\
{\bf (iv) The Frobenius semiradical}\\
Put \quad $F(L) = \sum\limits_{f \in L^\ast_{reg}} L(f)$.\quad This is a characteristic ideal of $L$ contained in $L_\Lambda$ and \quad $F(F(L)) = F(L)$.  We have the following characterization: \quad $Z(D(L)) \subset D(F(L))$\quad and if $L$ is algebraic then $F(L)$ is the smallest Lie subalgebra of $L$ with this property.  For example, if $L$ is square integrable then \quad $F(L) = Z(L)$\quad and so \quad $Z(D(L)) = D(Z(L))$. \quad In particular, \quad $Z(U(L)) = U(Z(L))$, \quad which is a polynomial algebra.  This phenomenon occurs quite frequently in section 5. Clearly, \quad $F(L) = 0$ \quad if and only if $L$ is Frobenius.  Therefore $F(L)$ is called the Frobenius semiradical of $L$.  At the other end of the spectrum we have the Lie algebras for which $F(L) = L$, which we call quasi quadratic.  These are unimodular and they do not possess any proper semi-invariants.  They form a large class, which include all quadratic Lie algebras (and hence all abelian and all semisimple Lie algebras) [O6].  This illustrates how far Frobenius and semisimple Lie algebras are removed from each other.\\
{\bf REMARK.}\quad Some time ago Mustapha Rais kindly sent us an unpublished manuscript by Andr\'e Cerezo [C2], in which the soul (respectively the rational soul) of a Lie algebra $L$ is introduced and studied.  This is the smallest Lie subalgebra of $L$ whose enveloping algebra (resp. enveloping quotient division ring) contains $Z(U(L))$ (resp. $Z(D(L)))$.  Clearly, the rational soul of $L$ coincides with $F(L)$ in case $L$ is algebraic.\\
\ \\
{\bf (v) The number \quad \boldmath{$c(L) = (\mbox{dim} L + i(L))/2$}}\\
$c(L)$, which is an integer, occurs frequently in the theory of Lie algebras.\\
\\
For instance:\\
{\bf THEOREM 2.5.} [OV] ($k$ algebraically closed)\\
Assume that $Sz(S(L))$ is freely generated by homogeneous elements \quad $f_1,\ldots, f_r$.  Then \quad $\sum\limits_{i=1}^r \deg f_i \leq c(L)$.\\
\ \\
The following, which is a minor extension of [OV, Corollary 1.2], is useful as a first test for coregularity.\\
\ \\
{\bf COROLLARY 2.6.} ($k$ algebraically closed)\\
Assume that $L$ has no proper semi-invariants in $S(L)$ and that \quad $S(L)^L = k[f_1,\ldots, f_r]$ \quad a polynomial algebra for some homogeneous $f_i \in S(L)$.  Then
$$3\ i(L) \leq \dim L + 2 \dim Z(L)$$
Moreover, if equality occurs then \quad $\sum\limits_{i=1}^r \deg f_i = c(L)$ \quad and \quad $\deg f_i \leq 2$, $i:1,\ldots, r$.\\
\ \\
On the other hand $c(L)$ is also the maximum transcendence degree of a Poisson commutative subfield of $R(L)$ [Sa].  See also [JL], [PY].\\
Similarly, $c(L)$ is an upperbound for the transcendence degree of maximal subfields of $D(L)$ and this bound can be achieved quite often.  However a  strictly smaller transcendence degree in $D(L)$ (and also in $R(L))$ is possible as the following example demonstrates.\\
Let $L$ be the solvable, almost algebraic Lie algebra over $k$ with basis \quad $\{x_0, x_1,\ldots x_n\}$, $n\geq 2$, with nonzero Lie brackets \quad $[x_0, x_i] = \lambda_ix_i$, $i:1,\ldots, n$, such that \quad $\lambda_1,\ldots,\lambda_n$ \quad are linearly independent over $\IQ$.  Clearly, $i(L) = n-1$\quad and \quad $c(L) = (n+1+n-1)/2=n$.  Then $k(x_0)$ is a maximal (resp. Poisson) commutative subfield of $D(L)$ (resp. of $R(L)$) [O5, Theorems 7, 12] and \quad $\mbox{tr}\deg_k (k(x_0)) = 1 < n = c(L)$.\\
\ \\
Next, let $L$ be an $n$-dimensional algebraic Lie algebra satisfying the Gelfand-Kirillov conjecture [GK], [BGR], [J1], [AOV1,2], [O7], i.e. $D(L)$ can be generated by elements \quad $z_1,\ldots, z_r$, $p_1,\ldots, p_s$, $q_1,\ldots, q_s \in D(L)$\quad such that \quad $[p_i,p_j] = 0$, $[q_i,q_j]=0$, $[p_i,q_j] = \delta_{ij}$, $i,j:1,\ldots, s$\quad and \quad $Z(D(L)) = k(z_1,\ldots, z_r)$, a purely transcendental extension of $k$. Then \quad $r=i(L)$, $n = r+2s$\quad and \quad $c(L) = r+s$\quad is the transcendence degree of the maximal commutative subfield \quad $k(z_1,\ldots, z_r, p_1,\ldots, p_s)$\quad of $D(L)$.\\
\ \\
Finally, let $P$ be a commutative Lie subalgebra of $L$ such that \quad $\dim P = c(L)$, i.e. $P$ is a commutative polarization (CP) w.r.t. any $f \in L^\ast_{reg}$.  Then $D(P)$ is a maximal subfield of $D(L)$ [O5, p. 706] with \quad $\mbox{tr} \deg_k(D(P)) = \dim P = c(L)$.\\
It is easy to see that \quad $F(L) \subset P$ \quad and hence is commutative.  $P$ plays a special role in the construction of the irreducible representations of $U(L)$ and their kernels, the primitive ideals.  If in addition, $P$ is an ideal of a completely solvable Lie algebra $L$, then $P$ is a Vergne polarization of $L$.  Furthermore, \quad $Sz(U(L)) \subset U(P)$\quad and \quad $P\subset L_\Lambda$ \quad [EO, p. 141].  In section 5 we will provide for each nilpotent Lie algebra a CP-ideal (if it exists).
\newpage
\begin{center}
{\bf 3. THE MAIN RESULT}
\end{center}
{\bf LEMMA 3.1.} Let $L$ be Frobenius and $H$ an ideal of $L$ such that $L/H$ is nilpotent.  Then
$$i(H) = \dim L - \dim H$$
{\bf PROOF.} Let \quad $f \in L^\ast_{reg}$ \quad such that its restriction \quad $f|_H\ \in H_{reg}^\ast$. Denote by $L^\infty$ the intersection of all terms $C^i(L)$ of the descending central series of $L$. With respect to the nondegenerate bilinear form \quad $B_f(x,y) = f([x,y])$, $x, y \in L$, we have that \quad $(L^{\infty})^{\perp} \subset L^\infty$ \quad [MO, Corollary 6.2].  Now, $L/H$ being nilpotent implies that \quad $L^\infty \subset H$.  Hence,
$$H^\perp \subset (L^\infty)^\perp \subset L^\infty \subset H$$
Consequently,
\begin{eqnarray*}
i(H) &=& \dim H(f|_H) = \dim (H \cap H^\perp) = \dim H^\perp\\
&=& \dim L - \dim H
\end{eqnarray*}
\hfill$\Box$\\
\\
{\bf PROPOSITION 3.2.} \quad Let $H$ be an ideal of $L$ such that $L/H$ is nilpotent.  Suppose $v \in S(L)$ is a semi-invariant with weight $\lambda \in L^\ast$.  Then $S(H)$ contains a semi-invariant $w$ of $S(L)$ with the same weight.\\
In particular, if $U(L)$ is primitive then each semi-invariant $v$ of $S(L)$ is contained in $S(H)$ and
$$Sz(S(L)) \subset Sz(S(H))$$

{\bf PROOF.} \quad Since $L/H$ is nilpotent, we can find ideals $H_i$ of $L$ with \quad $\dim H_i = i$, \quad $[L,H_i] \subset H_{i-1}$ \quad and such that
$$L = H_n \supset \ldots \supset H_i \supset H_{i-1} \supset \ldots \supset H_d = H$$
If \quad $v \in S(H)$\quad then put \quad $w = v$.  So, we may assume that \quad $v \in S(H_i)\backslash S(H_{i-1})$ \quad with $i > d$. Choose \quad $x \in H_i\backslash H_{i-1}$.  Then \quad $S(H_i) = S(H_{i-1}) [x]$.  So, $v$ can be written as a polynomial in $x$ with coefficients in $S(H_{i-1})$:
$$v = a_m x^m + \ldots + a_1 x + a_0,\ a_m \neq 0$$
Take any $y \in L$.  Then \quad $ad\ y (x) \in H_{i-1}$\quad and\\
$\lambda(y)v = ad\ y(v) = ad\ y(a_m)x^m+$ terms of lower degree in $x$.\\
\\
It follows that \quad $ad\ y(a_m) = \lambda(y)a_m$\quad for all $y \in L$, i.e. \quad $a_m \in S(L)_\lambda \cap S(H_{i-1})$.  After repeating the same reasoning a number of times, we find a semi-invariant $w$ of $S(L)$ with weight $\lambda \in \Lambda(L)$, contained in $S(H)$.\\
Now, suppose $U(L)$ is primitive, which is equivalent with \quad $R(L)^{L} = k$. [O3, p. 69].  Since $v$ and $w$ are semi-invariants with the same weight $\lambda \in \Lambda(L)$, we observe that 
$$ad\ y(vw^{-1}) = (w\ ad\ y(v) - v\ ad\ y(w))w^{-2} = 0$$
for all $y \in L$, i.e. \quad $vw^{-1} \in R(L)^L = k$ \quad Hence, for some \quad $a \in k$ \quad we have \quad $v = aw \in S(H)$. \hfill $\Box$\\
\ \\
{\bf LEMMA 3.3.} \quad Let \quad $\lambda \in L^\ast$ \quad be a weight.  Then
$$S(L)_\lambda = \bigoplus\limits_i (S^i (L) \cap S(L)_\lambda)$$
 
{\bf PROOF.} The inclusion $\supset$ is obvious.  For $\subset$ we take \quad $v \in S(L)_\lambda$ \quad and let 
$$v = v_p + \ldots + v_0$$
be its unique decomposition into homogeneous components (i.e. $v_i \in S^i(L))$.  Then for any $x \in L$: 
$$\sum\limits_{i=0}^p ad\ x (v_i) = ad\ x (v) = \lambda(x) v = \sum\limits_{i=0}^p \lambda(x) v_i$$
Since \quad $ad\ x (v_i) \in S^i(L)$\quad for all $i$ we may conclude that \quad $ad\ x(v_i) = \lambda(x) v_i$\quad for all $x \in L$, i.e. \quad $v_i \in S^i(L) \cap S(L)_\lambda$. \hfill $\Box$\\
\ \\
{\bf PROPOSITION 3.4.}\quad Let {\myfonts g} be a finite-dimensional Lie algebra over $k$.  Suppose there exists a split torus \quad $T \subset Der$ {\myfonts g}. Consider the semi-direct product \quad $L = T\ \oplus\ ${\myfonts g}.   Then,
$$Sz(S(\mbox{\myfonts g})) \subset Sz(S(L))\quad \mbox{and}\quad Sz(U(\mbox{\myfonts g}))\subset Sz(U(L))$$
Moreover, each weight \quad $\lambda \in \Lambda(\mbox{\myfonts g})$ \quad can be extended to a weight \quad $\lambda' \in \Lambda(L)$\quad i.e. $\lambda'|_{\mbox{\myfonts g}}\ = \lambda$.\\
\ \\
{\bf PROOF.} Take any weight $\lambda \in \Lambda({\mbox{\myfonts g}})$. By the previous lemma and the fact that $Sz(S(\mbox{{\myfonts g}})) = \bigoplus\limits_\lambda S({\mbox{\myfonts g}})_\lambda$, it suffices to show that for all $m$ \quad $V_{m,\lambda} = S^m ({\mbox{\myfonts g}}) \cap S({\mbox{\myfonts g}})_\lambda$ \quad is contained in $Sz(S(L))$.  First we notice that $S({\mbox{\myfonts g}})_\lambda$ is a $T$-submodule of $S({\mbox{\myfonts g}})$, since $S({\mbox{\myfonts g}})_\lambda$ is stable under the derivations of {\myfonts g} [Mon, p. 265].  Next, we choose a basis \quad $x_1,\ldots, x_n$\quad of {\myfonts g} consisting of common eigenvectors for all $t \in T$.  Then for each $m$ \quad $S^m({\mbox{\myfonts g}})$ \quad is a finite dimensional diagonalizable $T$-module (indeed the monomials \quad $x_1^{m_1} \ldots x_n^{m_n}$, \quad $\sum\limits_i m_i = m$, form a basis of $S^m ({\mbox{\myfonts g}})$ consisting of eigenvectors for any $t \in T$).  Hence, the same holds for the $T$-submodule $V_{m,\lambda}$.  Therefore if $V_{m,\lambda} \neq 0$ it contains a basis \quad $v_1,\ldots,v_p$ \quad such that \quad $t(v_i) = \mu_i (t) v_i$\quad for all \quad $t \in T$\quad for some \quad $\mu_i \in T^\ast$.  \quad On the other hand, \quad $ad\ x(v_i) = \lambda(x) v_i$\quad for all $x \in {\mbox{\myfonts g}}$.  Hence, each $v_i$ is a semi-invariant for $L = T\oplus {\mbox{\myfonts g}}$, say with weight $\lambda' \in \Lambda(L)$.  Note that $\lambda'$ is an extension of $\lambda$.  It follows that \quad $V_{m,\lambda} \subset Sz(S(L))$. \quad Consequently, \quad $Sz(S({\mbox{\myfonts g}})) \subset Sz(S(L))$,\quad which implies that \quad $Sz(U({\mbox{\myfonts g}})) \subset Sz(U(L))$\quad (take the image under the symmetrization $s$). \hfill $\Box$\\
\ \\
{\bf NOTATION.} Let $L$ be a Frobenius Lie algebra with basis \quad $x_1,\ldots, x_n$.\quad We recall that the Pfaffian \quad $Pf([x_i,x_j]) \in S(L)$\quad is a homogeneous semi-invariant of degree $\frac{1}{2}\dim L$ with weight \quad $\tau \in \Lambda(L)$\quad where \quad $\tau(x) = \mbox{tr}(ad\ x)$, $x \in L$.  Also,
$$(Pf([x_i,x_j]))^2 = \mbox{det} ([x_i,x_j])$$
We now put
$$Pf(L) = Pf([x_i, x_j])\quad \mbox{and}\quad \Delta(L) = \mbox{det}([x_i,x_j])$$
These are well-defined, up to a nonzero scalar multiple.\\
\ \\
{\bf THEOREM 3.5.} \quad Let {\myfonts g} be a finite dimensional Lie algebra over $k$.  Suppose there exists a split torus \quad $T \subset \mbox{Der}\ {\mbox{\myfonts g}}$ \quad such that the semi-direct product \quad $L = T \oplus {\mbox{\myfonts g}}$\quad is a Frobenius Lie algebra.  Let \quad $v_1,\ldots, v_r$ \quad be the irreducible factors of $\Delta(L)$.\\
Then the following hold:
\begin{itemize}
\item[1.] $Sz(S({\mbox{\myfonts g}})) = Sz(S(L)) = k[v_1, \ldots,v_r]$, a polynomial algebra.
\item[2.] $Sz(U({\mbox{\myfonts g}})) = Sz(U(L)) = k[s(v_1),\ldots, s(v_r)]$, a polynomial algebra.  This coincides with $Z(U(L_\Lambda))$ if $k$ is algebraically closed.
\item[3.] $\dim T = i({\mbox{\myfonts g}})$\quad and \quad $r= i(L_\Lambda) = \dim L - \dim L_\Lambda$
\item[4.] $\Lambda({\mbox{\myfonts g}}) = \{\lambda|_{\mbox{\myfonts g}}\ | \ \lambda \in \Lambda(L)\}$ \quad and \quad ${\mbox{\myfonts g}}_\Lambda = {\mbox{\myfonts g}}\cap L_\Lambda$
\item[5.] $\mbox{deg}\ v_1 + \ldots + \mbox{deg}\ v_r \leq c({\mbox{\myfonts g}}) = c(L)$
\item[6.] If $P$ is a CP of {\myfonts g}, then $P$ is also a CP of $L$.
\item[7.] If {\myfonts g} has no proper semi-invariants in $S({\mbox{\myfonts g}})$ (i.e. ${\mbox{\myfonts g}} = {\mbox{\myfonts g}}_\Lambda)$ then
$$S({\mbox{\myfonts g}})^{{\mbox{\myfonts g}}} = k[v_1,\ldots, v_r]\quad \mbox{and}\quad Z(U({\mbox{\myfonts g}})) = k[s(v_1),\ldots, s(v_r)]$$
If in addition $k$ is algebraically closed or if {\myfonts g} is algebraic, then \quad ${\mbox{\myfonts g}} = L_\Lambda$.
\end{itemize}

{\bf PROOF.}\quad (1) By the previous proposition \quad $Sz(S({\mbox{\myfonts g}})) \subset Sz(S(L))$\quad and each weight $\lambda \in \Lambda({\mbox{\myfonts g}})$\quad can be extended to a weight \quad $\lambda' \in \Lambda(L)$. On the other hand, $U(L)$ is primitive as $L$ is Frobenius.  Moreover, {\myfonts g} is an ideal of $L$ and $L/{\mbox{\myfonts g}}$ is abelian.  Then, by Proposition 3.2 each semi-invariant of $S(L)$ is already contained in $S({\mbox{\myfonts g}})$ and \quad $Sz(S(L)) \subset Sz(S({\mbox{\myfonts g}}))$.  Consequently, \quad $Sz(S({\mbox{\myfonts g}})) = Sz(S(L))$, the latter being a polynomial algebra in the variables \quad $v_1,\ldots, v_r$ \quad by Theorem 2.3.  Note also that if \quad $\mu \in \Lambda(L)$\quad then its restriction \quad $\mu|_{{\mbox{\myfonts g}}}\ \in \Lambda({\mbox{\myfonts g}})$.\\
\\
(2) This follows immediately from (1) and the fact that \quad $s : Sz(S(L)) \rightarrow Sz(U(L))$\quad is an algebra isomorphism which maps $Sz(S({\mbox{\myfonts g}}))$ onto $Sz(U({\mbox{\myfonts g}}))$.\\
\\
(3) Both quotients $L/{\mbox{\myfonts g}}$ and $L/L_\Lambda$ [DNO, p. 330] are abelian.  Hence, by Lemma 3.1 and Theorem 2.3:
$$i({\mbox{\myfonts g}}) = \dim L - \dim {\mbox{\myfonts g}} = \dim T\quad \mbox{and}\quad i(L_\Lambda) = \dim L - \dim L_\Lambda = r$$

(4) From the proof of (1) we deduce that the map \quad $\Lambda(L) \rightarrow \Lambda({\mbox{\myfonts g}})$, sending $\lambda$ into its restriction $\lambda|_{\mbox{\myfonts g}}$, is surjective.  It follows that
$$\Lambda({\mbox{\myfonts g}}) = \{\lambda|_{\mbox{\myfonts g}}\ | \ \lambda \in \Lambda(L)\}\quad\mbox{and}\quad {\mbox{\myfonts g}}_\Lambda = {\mbox{\myfonts g}} \cap L_{\Lambda}$$

(5) $v_1,\ldots, v_r$ \quad are also the irreducible factors of $Pf(L)$ as $\quad \Delta(L)=(Pf(L))^2$.  Hence, the sum of their degrees is at most:
\begin{eqnarray*}
\mbox{deg}(Pf(L)) &=& \frac{1}{2}\dim L = \frac{1}{2} (\dim {\mbox{\myfonts g}} + \dim T)\\
&=& \frac{1}{2}(\dim {\mbox{\myfonts g}} + i({\mbox{\myfonts g}})) = c({\mbox{\myfonts g}})
\end{eqnarray*}
On the other hand, $c(L) = \frac{1}{2} (\dim L + i(L)) = \frac{1}{2} \dim L = c({\mbox{\myfonts g}})$\\
\\
(6) Let $P$ be a CP of {\myfonts g}, i.e. $P$ is a commutative Lie subalgebra of {\myfonts g} for which $\dim P = c({\mbox{\myfonts g}})$.  Clearly, $P$ is also a CP of $L$ since \quad $c({\mbox{\myfonts g}}) = c(L)$.\\
\\
(7) $S({\mbox{\myfonts g}})^{{\mbox{\myfonts g}}} = Sz(S({\mbox{\myfonts g}}))$\quad and \quad $Z(U({\mbox{\myfonts g}})) = Sz(U({\mbox{\myfonts g}}))$\quad as {\myfonts g} does not have proper semi-invariants in $S({\mbox{\myfonts g}})$ nor in $U({\mbox{\myfonts g}})$.  So, it suffices to apply (1) and (2).\\
Next, we assume in addition that either $k$ is algebraically closed or that {\myfonts g} is algebraic.  Using (4), we deduce from \quad ${\mbox{\myfonts g}} = {\mbox{\myfonts g}}_\Lambda$ \quad that \quad ${\mbox{\myfonts g}} \subset L_\Lambda$.  Hence it is enough to show that \quad $\dim {\mbox{\myfonts g}} = \dim L_\Lambda$.  As there are no proper semi-invariants in $U({\mbox{\myfonts g}})$ we see that $Z(D({\mbox{\myfonts g}}))$ is the quotient field of $Z(U({\mbox{\myfonts g}}))$.  Therefore,
$$Z(D({\mbox{\myfonts g}})) = k(s(v_1),\ldots, s(v_r))$$
In particular, \quad $r = \mbox{tr deg}_k Z(D({\mbox{\myfonts g}})) = i({\mbox{\myfonts g}})$ by Theorem 2.1.  So,
$$\dim L - \dim L_\Lambda = \dim L - \dim {\mbox{\myfonts g}}$$
which yields \quad $\dim L_\Lambda = \dim {\mbox{\myfonts g}}$. \hfill $\Box$\\
\ \\
{\bf REMARK.} (1) and (2) of the theorem remain valid under the condition that $U(L)$ is primitive (instead of requiring $L$ to be Frobenius).  In that case \quad $v_1,\ldots, v_r$ \quad are the irreducible factors of $\Delta_n(L)$, a special semi-invariant of $S(L)$ [DNO, p. 337].\\
\ \\
{\bf COROLLARY 3.6.} \quad See also [J2].  Let $k$ be algebraically closed and let $L$ be a simple Lie algebra over $k$ of rank $r$ with triangular decomposition \quad $L = N^- \oplus H \oplus N$.  Consider the Borel subalgebra \quad $B = H \oplus N$.  If $L$ is not of one of the following types $A_n$, $n \geq 2$; $D_{2t+1}$, $t \geq 2$; $E_6$, then
\begin{itemize}
\item[(1)] $B$ is Frobenius
\item[(2)] $i(N) = r$, \quad $c(B) = c(N)$, \quad $N = B_\Lambda$
\item[(3)] $Sz(U(B)) = Z(U(N))$, a polynomial algebra in $r$ generators which can be explicitly determined by the method of Theorem 3.5.
\end{itemize}
{\bf PROOF.} \quad Clearly, $B$ can be considered as the semi-direct product of the split torus \quad $T = ad_NH \subset \mbox{Der} N$ \quad with the nilradical $N$.  Since $B$ is Frobenius [EO,~p. 146] the previous theorem can be applied. \hfill $\Box$
\newpage
\begin{center}
{\bf 4. EXAMPLES}
\end{center}
1. Let {\myfonts g} be the 5-dimensional solvable Lie algebra over $k$ with basis $x_1,\ldots, x_5$ and with nonvanishing Lie brackets: $[x_1, x_3] = x_3$ \quad $[x_1, x_4] = x_4$\quad $[x_1, x_5] = x_5$ \quad $[x_2, x_3] = x_4$ \quad $[x_2, x_4] = x_5$.\\
We want to apply Theorem 3.5.  Since $i({\mbox{\myfonts g}}) = 1$ we need to find a 1-dimensional split torus \quad $T \subset Der\ {\mbox{\myfonts g}}$ such that the semi-direct product \quad $L = T \oplus {\mbox{\myfonts g}}$ \quad is Frobenius.  Put \quad $T = <t>$ \quad with \quad $t = \diag (0,-1,2,1,0)$.  The matrix of Lie brackets of $L$ is:
\begin{eqnarray*}
\begin{array}{c|cccccc}
 &t &x_1 &x_2 &x_3 &x_4 &x_5\\
\hline
t &0 &0 &-x_2 &2x_3 &x_4 &0\\
x_1 &0 &0 &0 &x_3 &x_4 &x_5\\
x_2 &x_2 &0 &0 &x_4 &x_5 &0\\
x_3 &-2x_3 &-x_3 &-x_4 &0 &0 &0\\
x_4 &-x_4 &-x_4 &-x_5 &0 &0 &0\\
x_5 &0 &-x_5 &0 &0 &0 &0
\end{array}
\end{eqnarray*}
Its determinant is \quad $\Delta(L) = x_5^2 (x_4^2 - 2x_3 x_5)^2$\quad which is nonzero.  So, $L = T\oplus {\mbox{\myfonts g}}$ is Frobenius.  Clearly, $x_5$ and $x_4^2 - 2x_3 x_5$ are the irreducible factors of $\Delta(L)$.  Hence, by Theorem 3.5:
$$Sz(S({\mbox{\myfonts g}})) = Sz(S(L)) = k[x_5, x_4^2 - 2x_3x_5]$$
which is a polynomial algebra.\\
Finally, by symmetrization:
$$Sz(U({\mbox{\myfonts g}})) = Sz(U(L)) = k[x_5, x_4^2-2x_3x_5]$$
2. The following nonsolvable examples, selected from [AOV2], have a Levi decomposition \quad ${\mbox{\myfonts g}} = S\oplus R$, where \quad $S = sl(2,k)$, with standard basis \quad $h,x,y$\quad and nonzero Lie brackets \quad $[h,x] = 2x$, \quad $[h,y] = -2y$,\quad $[x,y] = h$.\\
$e_0, e_1,\ldots, e_p$ \quad will be a basis of $R$. $W_n$ will be the $(n+1)$-dimensional irreducible $sl(2,k)$-module, with standard basis $e_0, e_1, \ldots, e_n$.\\
\ \\
(i) ${\mbox{\myfonts g}} = sl(2,k) \oplus W_1$, with $W_1$ abelian [AOV2, p. 554]. {\myfonts g} is a 5-dimensional algebraic Lie algebra of index 1.  Put \quad $T = <t>$, \quad $t = \diag(0,0,0,1,1)$\quad and \quad $L = T \oplus {\mbox{\myfonts g}}$.\\
Its matrix of Lie brackets
\begin{eqnarray*}
\begin{array}{c|cccccc}
 &t &h &x &y &e_0 &e_1\\
\hline
t &0 &0 &0 &0 &e_0 &e_1\\
h &0 &0 &2x &-2y &e_0 &-e_1\\
x &0 &-2x &0 &h &0 &e_0\\
y &0 &2y &-h &0 &e_1 &0\\
e_0 &-e_0 &-e_0 &0 &-e_1 &0 &0\\
e_1 &-e_1 &e_1 &-e_0 &0 &0 &0
\end{array}
\end{eqnarray*}
has determinant \quad $\Delta(L) = 4(e_1^2 x + e_0e_1h - e_0^2 y)^2 \neq 0$.\\
So, $L$ is Frobenius.  On the other hand, {\myfonts g} has no proper semi-invariants in $S({\mbox{\myfonts g}})$ as \quad $[{\mbox{\myfonts g}},{\mbox{\myfonts g}}]={\mbox{\myfonts g}}$. By (7) of Theorem 3.5 we conclude that \quad ${\mbox{\myfonts g}} = L_\Lambda$ \quad and 
$$Z(U({\mbox{\myfonts g}})) = Sz(U(L)) = k[e_1^2x + e_0 e_1 h - e_0^2y]$$
(ii) ${\mbox{\myfonts g}} = sl(2,k) \oplus R$, where \quad $R = W_2 \oplus W_1$\quad with standard basis \quad $e_0, e_1, e_2; e_3, e_4$. [AOV2, p. 567], [O6, p. 908].  The nontrivial action of $sl(2,k)$ on $R$ is given by\\
\ \\
$[h, e_0] = 2e_0\quad [h,e_2] = -2e_2 \quad \quad [h, e_3] = e_3 \quad [h, e_4] = -e_4$\\
$[x, e_1] = 2e_0\quad [x, e_2] = e_1 \quad [x, e_4] = e_3$ $[y, e_0] = e_1 \quad [y, e_1] = 2e_2 \quad [y, e_3] = e_4$.\\
{\mbox{\myfonts g}} is an $8$-dimensional algebraic Lie algebra of index 2.  Consider the split torus \quad $T = <t_1, t_2>$; \quad $t_1 = \diag (0,0,0,1,1,1,0,0)$, \quad $t_2 = \diag (0,0,0,0,0,0,1,1)$ \quad and \quad $L = T \oplus {\mbox{\myfonts g}}$.\\
Then \quad $\Delta(L) = 4(e_1^2 - 4e_0e_2)^2 (e_0e_4^2 - e_1e_3e_4 + e_2 e_3^2)^2 \neq 0$.\
Hence, $L$ is a 10-dimensional Frobenius Lie algebra.  Clearly, $[{\mbox{\myfonts g}},{\mbox{\myfonts g}}] = {\mbox{\myfonts g}}$.  By (7) of Theorem 3.5 we obtain that \quad ${\mbox{\myfonts g}} = L_\Lambda$ \quad and
$$Z(U({\mbox{\myfonts g}})) = Sz(U(L)) = k[e_1^2 - 4e_0e_2, e_0e_4^2 - e_1 e_3e_4 + e_2e_3^2]$$
(iii) Let $R$ be the 5-dimensional nilpotent Lie algebra with basis \quad $e_0, e_1, e_2, e_3, e_4$ \quad and nonzero Lie brackets $[e_2, e_4] = e_0$ \quad $[e_3, e_4] = e_1$. \quad $sl(2,k)$ acts on $R$ as follows: $[h, e_0] = e_0$ \quad $[h, e_1] = -e_1$ \quad $[h,e_2] = e_2$ \quad $[h, e_3] = -e_3$ \quad $[x, e_1] = e_0$ \quad $[x, e_3] = e_2$ \quad $[y,e_0] = e_1$ \quad $[y, e_2] = e_3$.\\
Consider the semi-direct product \quad ${\mbox{\myfonts g}} = sl(2,k) \oplus R$. {\mbox{\myfonts g}} is an 8-dimensional algebraic Lie algebra of index 2 [AOV2, p. 573].\\
Put \quad $T = <t_1, t_2>$;\quad $t_1 = diag (0,0,0,1,1,0,0,1)$, $t_2 = diag (0,0,0,1,1,1,1,0)$.  Then,
$$\Delta(L) = 4(e_0e_3-e_1e_2)^2 (e_1^2 x + e_0e_1h - e_0^2y - e_0e_3 e_4 + e_1 e_2 e_4)^2 \neq 0$$
Hence, $L$ is a 10-dimensional Frobenius Lie algebra.  Clearly, {\mbox{\myfonts g}} has no proper semi-invariants in $S({\mbox{\myfonts g}})$ as $R$ is nilpotent.  By (7) of Theorem 3.5 we may conclude that \quad ${\mbox{\myfonts g}} = L_\Lambda$ \quad and after symmetrization
$$Z(U({\mbox{\myfonts g}})) = Sz(U(L)) = k[e_0e_3 - e_1e_2, e_1^2x + e_0e_1h - e_0^2y - e_0e_3e_4 + e_1e_2e_4]$$

{\bf REMARK.}\quad In the last 3 examples we notice that the algebraically independent generators of $Z(D({\mbox{\myfonts g}}))$ we obtained in [AOV2] and [O6] (needed for the verification of the GK-conjecture) turn out to be the generators of $Z(U({\mbox{\myfonts g}}))$.

\begin{center}
{\bf 5. INDECOMPOSABLE NILPOTENT LIE ALGEBRAS OF DIMENSION \boldmath{$\leq 7$ $(k = \IC)$}}
\end{center}
The main purpose is to describe \quad $Z = Z(U({\mbox{\myfonts g}}))$ \quad for each Lie algebra {\mbox{\myfonts g}}, but also to give the Frobenius semiradical \quad $F = F({\mbox{\myfonts g}})$ \quad and if they exist a CP-ideal (CPI) and a torus \quad $T \subset \mbox{Der}\ {\mbox{\myfonts g}}$ \quad for which the semi-direct product \quad $T \ \oplus \ {\mbox{\myfonts g}}$ \quad is Frobenius.  Sometimes the possible solutions for $T$ form a family depending on parameters $\in k$.  Other abbreviations are: $i = i({\mbox{\myfonts g}})$, $r = \mbox{rank}\ {\mbox{\myfonts g}}$, $c = c({\mbox{\myfonts g}})$, SQ.I. = square integrable, $Q(Z)=$ the quotient field of $Z(U({\mbox{\myfonts g}}))$.  For dimension $\leq 6$ all Lie algebras are listed, while in dimension 7 only those for which $i({\mbox{\myfonts g}})\leq  \mbox{rank}\ {\mbox{\myfonts g}}$.\\
These Lie algebras are coregular (i.e. $Z(U({\mbox{\myfonts g}}))$ is a polynomial algebra over $k$), except ${\mbox{\myfonts g}}_{5,5}$, ${\mbox{\myfonts g}}_{6,6}$, ${\mbox{\myfonts g}}_{6,16}$ and ${\mbox{\myfonts g}}_{7,3.17}$.\\
\ \\
{\bf Notation:} We use the same notation as Magnin [Ma2] and Carles [Ca].  In addition we include the notation of Morozov (M$\_$) [Mor] and Cerezo (C$\_$) [C2] in dimension 6 and of Seeley ((...)) [Se] and Romdhani (R$\_$)[R] in dimension 7.\\
\ \\
{\bf \boldmath{$\dim {\mbox{\myfonts g}} \leq 5$}} (See also [D3, pp. 322-330])
\begin{itemize}
\item[1.] $\mbox{\myfonts g}_3$ (3-dim Heisenberg Lie algebra)\\
$[x_1, x_2] = x_3$.\\
SQ.I.\quad $i = 1\quad r = 2 \quad c = 2 \quad Z= k[x_3] \quad F=<x_3>$\\
$CPI = <x_2, x_3>$\\
$T = <t>$, \quad $t = \diag (\alpha, 1-\alpha, 1).$
\item[2.] $\mbox{\myfonts g}_4$ (4-dim standard filiform Lie algebra)\\
$[x_1, x_2] = x_3 \quad [x_1, x_3] = x_4$.\\
$i = 2\quad r = 2 \quad Z= k[x_4, x_3^2 - 2x_2x_4]\\
F=<x_2,x_3,x_4> = CPI$\\
$T = <t_1,t_2>$, \quad $t_1 = \diag(0,1,1,1)$, \quad $t_2 = \diag (1,-2,-1,0)$.\\
\end{itemize}

{\bf \boldmath{$\dim {\mbox{\myfonts g}} = 5$, $i(\mbox{\myfonts g}) \leq rank \ \mbox{\myfonts g}$}}

\begin{itemize}
\item[3.] $\mbox{\myfonts g}_{5,1}$ (5-dim. Heisenberg Lie algebra)\\
$[x_1,x_3] = x_5 \quad [x_2, x_4] = x_5$.\\
SQ.I. $i = 1\quad r= 3 \quad c=3 \quad Z= k[x_5]$\\
$F = <x_5> \quad CPI = <x_3, x_4, x_5>$\\
$T = <t>, \quad t = \diag (\alpha,\beta, 1-\alpha, 1-\beta, 1)$.

\item[4.] $\mbox{\myfonts g}_{5,2}$\\
$[x_1,x_2] = x_4 \quad [x_1, x_3] = x_5$.\\
$i=3 \quad r=3 \quad c= 4 \quad Z = k[x_4, x_5, x_2x_5 - x_3x_4]$\\
$F = <x_2, x_3, x_4, x_5> = CPI$\\
$T = <t_1, t_2, t_3>$, \quad $t_1 = \diag(1,0,0,1,1),\\
t_2 = \diag (0,1,0,1,0), \quad t_3 = \diag (1, -1, -1, 0,0)$.

\item[5.] $\mbox{\myfonts g}_{5,3}$\\
$[x_1,x_2] = x_4 \quad [x_1,x_4] = x_5 \quad [x_2, x_3] = x_5$\\
SQ.I. $i = 1 \quad r= 2 \quad c = 3 \quad Z = k[x_5] \quad F = <x_5>\quad CPI = <x_3, x_4, x_5>$\\
$T = <t>$, \quad $t = \diag (\alpha, 1-2\alpha, 2\alpha, 1-\alpha, 1)$.

\item[6.] $\mbox{\myfonts g}_{5,6}$\\
$[x_1,x_2] = x_3 \quad [x_1, x_3] = x_4 \quad [x_1, x_4] = x_5 \quad [x_2, x_3] = x_5$\\
SQ.I. $i=1\quad r=1\quad c=3 \quad Z = k[x_5] \quad F = <x_5> \quad CPI = <x_3, x_4, x_5>$\\
$T = <t>$, \quad $t = \diag(1,2,3,4,5)$.
\end{itemize}

{\boldmath{$\dim \mbox{\myfonts g} =5, \quad i(\mbox{\myfonts g}) > \mbox{rank}\ \mbox{\myfonts g}$}}

\begin{itemize}
\item[7.] $\mbox{\myfonts g}_{5,4}$\\
$[x_1,x_2] = x_3 \quad [x_1, x_3] = x_4 \quad [x_2, x_3] = x_5$.\\
$i = 3 \quad r=2 \quad c=4 \quad Z = k[x_4, x_5, x_3^2 + 2x_1x_5 - 2x_2x_4]$\\
$F= \mbox{\myfonts g}_{5,4}$ \quad No CP's.

\item[8.] $\mbox{\myfonts g}_{5,5}$ (5-dim. standard filiform Lie algebra)\\
$[x_1, x_2] = x_3 \quad [x_1, x_3] = x_4 \quad [x_1, x_4] = x_5$.\\
$i = 3 \quad r = 2 \quad c = 4 \quad F = <x_2,x_3,x_4,x_5> = CPI$\\
$Z = k[x_5, f_1, f_2, f_3]\quad f_1 = 2x_3x_5 - x_4^2\quad f_2 = 3x_2x_5^2 - 3x_3 x_4 x_5 + x_4^3$ \\
$f_3 = 9x_2^2x_5^2 - 18x_2x_3x_4x_5 + 6x_2x_4^3 + 8x_3^3x_5 - 3x_3^2x_4^2$\\
Relation: \quad $f_1^3 + f_2^2 - x_5^2f_3 = 0$ \quad $Q(Z) = k(x_5, f_1, f_2)$.
\end{itemize}

{\boldmath{$\dim \mbox{\myfonts g} = 6, \quad i(\mbox{\myfonts g}) \leq \mbox{rank}\ \mbox{\myfonts g}$}}

\begin{itemize}
\item[9.] $\mbox{\myfonts g}_{6,1} \cong M_4 \cong C_{24}$\\
$[x_1,x_2] = x_5 \quad [x_1, x_4] = x_6 \quad [x_2, x_3] = x_6$.\\
SQ.I. $i = 2 \quad r = 3\quad c=4 \quad Z = k[x_5, x_6]$\\
$F = <x_5, x_6>  \quad CPI = <x_3, x_4, x_5, x_6>$\\
$T = <t_1,t_2>$, \quad $t_1 = \diag(\alpha, 1-\alpha, \alpha, 1-\alpha, 1,1), \quad t_2 = \diag (1,0,0,-1,1,0)$.

\item[10.] $\mbox{\myfonts g}_{6,2} \cong M_{12} \cong C_{22}$\\
$[x_1,x_2] = x_5 \quad [x_1, x_5] = x_6 \quad [x_3, x_4] = x_6$.\\
$i = 2 \quad r=3 \quad c=4\quad Z = k[x_6, x_5^2 - 2x_2x_6]$\\
$F=<x_2, x_5, x_6> \quad CPI = <x_2, x_4, x_5, x_6>$\\
$T = <t_1,t_2>$, \quad $t_1 = \diag(0,1,\alpha,1-\alpha, 1,1)$, \quad $t_2 = \diag(1,0,2,0,1,2)$.

\item[11.] $\mbox{\myfonts g}_{6,4} \cong M_7\cong C_{18}$\\
$[x_1,x_2] = x_4 \quad [x_1, x_3] = x_6\quad [x_2, x_4] = x_5$.\\
SQ.I. $i = 2 \quad r = 3 \quad c = 4\quad Z = k[x_5, x_6]$\\
$F = <x_5,x_6> \quad CPI = <x_3, x_4, x_5, x_6>$\\
$T = <t_1,t_2>, t_1 = \diag (1-2\alpha, \alpha, 2\alpha, 1-\alpha, 1,1)$, \quad $t_2 = \diag (0,1,0,1,2,0)$.

\item[12.] $\mbox{\myfonts g}_{6,5} \cong M_8 \cong C_{12}$\\
$[x_1,x_2] = x_4 \quad [x_1, x_4] = x_5 \quad [x_2, x_3] = x_6 \quad [x_2, x_4] = x_6$.\\
SQ.I. \quad $i = 2 \quad r=2 \quad c = 4 \quad Z = k[x_5, x_6]$\\
$F = <x_5, x_6>\quad CPI = <x_3, x_4, x_5, x_6>$\\
$T = <t_1, t_2>, \quad t_1 = \diag (1,0,1,1,2,1), \quad t_2 = \diag (2,-1,1,1,3,0)$.

\item[13.] $\mbox{\myfonts g}_{6,7} \cong M_{6} \cong C_{19}$\\
$[x_1, x_2] = x_4 \quad [x_1, x_3] = x_5 \quad [x_1, x_4] = x_6\quad [x_2, x_3] = -x_6$.\\
SQ.I. \quad $i = 2 \quad r=2 \quad c = 4 \quad Z = k[x_5, x_6]$\\
$F = <x_5,x_6>\quad CPI = <x_3, x_4, x_5, x_6>$\\
$T = <t_1, t_2>, \quad t_1 = \diag (0,1,0,1,0,1), \quad t_2 = \diag(1,-2,2,-1,3,0)$.

\item[14.] $\mbox{\myfonts g}_{6,8} \cong M_{9} \cong C_{13}$\\
$[x_1, x_2] = x_4 \quad [x_1, x_4] = x_5 \quad [x_2, x_3] = x_5\quad [x_2, x_4] = x_6$.\\
SQ.I. \quad $i = 2 \quad r=2 \quad c = 4 \quad Z = k[x_5, x_6]$\\
$F = <x_5,x_6>\quad CPI = <x_3, x_4, x_5, x_6>$\\
$T = <t_1, t_2>, \quad t_1 = \diag (1,0,2,1,2,1), \quad t_2 = \diag(2,-1,4,1,3,0)$.

\item[15.] $\mbox{\myfonts g}_{6,9} \cong M_{14} \cong C_{16}$\\
$[x_1, x_2] = x_4 \quad [x_1, x_3] = x_5 \quad [x_2, x_5] = x_6\quad [x_3, x_4] = x_6$.\\
$i = 2 \quad r=3 \quad c = 4 \quad Z = k[x_6, x_1x_6 + x_4x_5]$\\
$F = <x_1, x_4, x_5, x_6> = CPI$\\
$T = <t_1, t_2>, \quad t_1 = \diag (1,\alpha, -\alpha, 1+ \alpha,1-\alpha, 1), \quad t_2 = \diag(1,0,-1,1,0,0)$.

\item[16.] $\mbox{\myfonts g}_{6,10} \cong M_{13} \cong C_{17}$\\
$[x_1, x_2] = x_4 \quad [x_1, x_3] = x_5 \quad [x_1, x_4] = x_6\quad [x_3, x_5] = x_6$.\\
$i = 2 \quad r=2 \quad c = 4 \quad Z = k[x_6, x_4^2-2x_2x_6]$\\
$F = <x_2,x_4,x_6>\quad CPI = <x_2, x_4, x_5, x_6>$\\
$T = <t_1, t_2>, \quad t_1 = \diag (1,1,1,2,2,3), \quad t_2 = \diag(2,-4,-1,-2,1,0)$.

\item[17.] $\mbox{\myfonts g}_{6,11} \cong M_{17} \cong C_{11}$\\
$[x_1, x_2] = x_4 \quad [x_1, x_4] = x_5 \quad [x_1, x_5] = x_6\quad [x_2, x_3] = x_6$.\\
$i = 2 \quad r=2 \quad c = 4 \quad Z = k[x_6, x_5^2 -2x_4x_6]$\\
$F = <x_4,x_5,x_6>\quad CPI = <x_3, x_4, x_5, x_6>$\\
$T = <t_1, t_2>, \quad t_1 = \diag (1,0,3,1,2,3), \quad t_2 = \diag(-1,3,-3,2,1,0)$.

\item[18.] $\mbox{\myfonts g}_{6,13} \cong M_{16} \cong C_{9}$\\
$[x_1, x_2] = x_4 \quad [x_1, x_4] = x_5 \quad [x_1, x_5] = x_6\quad [x_2, x_3] = x_5\quad [x_3, x_4] = -x_6$.\\
$i = 2 \quad r=2 \quad c = 4 \quad Z = k[x_6, x_5^3 - 3x_4x_5x_6+ 3x_2x_6^2]$\\
$F = <x_2, x_4,x_5,x_6> = CPI$\\
$T = <t_1, t_2>, \quad t_1 = \diag (1,-2,2,-1,0,1), \quad t_2 = \diag(1,-3,2,-2,-1,0)$.

\item[19.] $\mbox{\myfonts g}_{6,14} \cong M_{11} \cong C_{8}$\\
$[x_1, x_2] = x_3 \quad [x_1, x_3] = x_4 \quad [x_1, x_4] = x_5\quad [x_2, x_3] = x_6$.\\
SQ.I. \quad $i = 2 \quad r=2 \quad c = 4 \quad Z = k[x_5, x_6]$\\
$F = <x_5,x_6>\quad CPI = <x_3, x_4, x_5, x_6>$\\
$T = <t_1, t_2>, \quad t_1 = \diag (1,0,1,2,3,1), \quad t_2 = \diag(2,-1,1,3,5,0)$.

\item[20.] $\mbox{\myfonts g}_{6,15} \cong M_{18}(-1) \cong C_{7}$\\
$[x_1, x_2] = x_3 \quad [x_1, x_3] = x_4 \quad [x_1, x_5] = x_6\quad [x_2, x_3] = x_5 \quad [x_2,x_4] = x_6$.\\
$i = 2 \quad r=2 \quad c = 4 \quad Z = k[x_6, x_3x_6 - x_4 x_5]$\\
$F = <x_3, x_4,x_5,x_6>= CPI $

\item[21.] $\mbox{\myfonts g}_{6,18} \cong M_{21} \cong C_{2}$\\
$[x_1, x_2] = x_3 \quad [x_1, x_3] = x_4 \quad [x_1, x_4] = x_5\quad [x_2, x_5] = x_6 \quad [x_3, x_4] = -x_6$.\\
$i = 2 \quad r=2 \quad c = 4 \quad Z = k[x_6, x_4^2 -2x_1x_6 - 2x_3x_5]$\\
$F = <x_1,x_3,x_4,x_5,x_6>$\quad No CP's\\
$T = <t_1, t_2>, \quad t_1 = \diag (1,-1,0,1,2,1), \quad t_2 = \diag(2,-3,-1,1,3,0)$.
\end{itemize}

{\boldmath{$\dim \mbox{\myfonts g} = 6,\quad i(\mbox{\myfonts g}) > \mbox{rank}\ \mbox{\myfonts g}$}}

\begin{itemize}
\item[22.] $\mbox{\myfonts g}_{6,3} \cong M_{3} \cong C_{21}$\\
$[x_1, x_2] = x_4 \quad [x_1, x_3] = x_5 \quad [x_2, x_3] = x_6$.\\
$i = 4 \quad r=3 \quad c = 5 \quad Z = k[x_4,x_5, x_6, x_1x_6-x_2x_5+x_3x_4]$\\
$F = \mbox{\myfonts g}_{6,3}$ \quad No CP's.

\item[23.] $\mbox{\myfonts g}_{6,6} \cong M_{1}\cong C_{20}$\\
$[x_1, x_2] = x_4 \quad [x_2, x_3] = x_6 \quad [x_2, x_4] = x_5$.\\
$i = 4 \quad r=3 \quad c = 5 \quad F = <x_1, x_3,x_4, x_5,x_6>= CPI $\\
$Z = k[x_5, x_6, f_1, f_2, f_3], \quad f_1 = x_4^2 + 2x_1x_5, \quad f_2 = x_3x_5 - x_4x_6,$\\
$f_3 = 2x_1x_6^2 + 2x_3x_4x_6 - x_3^2x_5$,\\
Relation: \quad $x_6^2f_1 - f_2^2 - x_5f_3 = 0\quad Q(Z) = k(x_5,x_6, f_1, f_2)$.

\item[24.] $\mbox{\myfonts g}_{6,12} \cong M_{15} \cong C_{10}$\\
$[x_1, x_2] = x_4 \quad [x_1, x_4] = x_5\quad [x_1, x_5] = x_6\quad [x_2, x_3] = x_6 \quad [x_2,x_4] = x_6$.\\
$i = 2 \quad r=1 \quad c = 4 \quad Z = k[x_6, x_5^2 + 2x_3x_6 - 2x_4x_6]$\\
$F = <x_3-x_4, x_5,x_6> \quad CPI = <x_3,x_4,x_5,x_6>$.

\item[25.] $\mbox{\myfonts g}_{6,16} \cong M_{2}\cong C_{5}$ (6-dim standard filiform Lie algebra)\\
$[x_1, x_2] = x_3 \quad [x_1, x_3] = x_4 \quad [x_1, x_4] = x_5 \quad [x_1, x_5] = x_6$.\\
$i = 4 \quad r=2 \quad c = 5 \quad F = <x_2, x_3,x_4, x_5,x_6>= CPI $\\
$Z = k[x_6, f_1, f_2, f_3, f_4],\\
f_1 = x_5^2 - 2x_4x_6,\\
f_2 = x_5^3 - 3x_4x_5x_6 + 3x_3x_6^2,\\
f_3 = x_4^2 + 2x_2x_6 - 2x_3x_5,\\
f_4 = 2x_4^3 + 6x_2x_5^2 + 9x_3^2x_6 - 12x_2x_4x_6 - 6x_3x_4x_5$.\\
Relation: \quad $f_1^3 - f_2^2 - 3x_6^2f_1f_3 + x_6^3f_4 = 0\\
Q(Z) = k(x_6, f_1, f_2,f_3)$.

\item[26.] $\mbox{\myfonts g}_{6,17} \cong M_{19} \cong C_{4}$\\
$[x_1, x_2] = x_3 \quad [x_1, x_3] = x_4 \quad [x_1, x_4] = x_5\quad [x_1, x_5] = x_6 \quad [x_2,x_3] = x_6.\\
i = 2 \quad r=1 \quad c = 4\quad Z = k[x_6, x_5^2 - 2x_4x_6]$\\
$F = <x_4,x_5,x_6>$\quad $CPI = <x_3, x_4, x_5, x_6>$

\item[27.] $\mbox{\myfonts g}_{6,19} \cong M_{20} \cong C_{3}$\\
$[x_1, x_2] = x_3 \quad [x_1, x_3] = x_4 \quad [x_1, x_4] = x_5\quad [x_1, x_5] = x_6$\\
$[x_2, x_3] = x_5 \ [x_2, x_4] = x_6$.\\
$i = 2 \quad r=1 \quad c = 4\quad Z = k[x_6, x_5^3 -3x_4x_5x_6 + 3x_3x_6^2]\\
F = <x_3, x_4, x_5,x_6> = CPI.$

\item[28.] $\mbox{\myfonts g}_{6,20} \cong M_{22} \cong C_{1}$\\
$[x_1, x_2] = x_3 \quad [x_1, x_3] = x_4 \quad [x_1, x_4] = x_5\quad [x_2, x_3] = x_5 \quad [x_2, x_5] = x_6$\\
$[x_3, x_4] = -x_6$.\\
$i = 2 \quad r=1 \quad c = 4\quad Z = k[x_6, 2x_5^3 + 3x_4^2 x_6 -6x_3x_5x_6 - 6x_1x_6^2]\\
F = <x_1, x_3, x_4, x_5,x_6>$ No CP's.
\end{itemize}

{\boldmath{$\dim \mbox{\myfonts g} = 7, \quad i(\mbox{\myfonts g}) \leq \mbox{rank} \ \mbox{\myfonts g}$}}

\begin{itemize}
\item[29.] $\mbox{\myfonts g}_{7,1.03} \cong (13457G) \cong R_{26}$\\
$[x_1, x_2] = x_3 \quad [x_1, x_3] = x_4 \quad [x_1, x_4] = x_5\quad [x_1, x_6] = x_7 \quad [x_2, x_3] = x_6$\\
$[x_2, x_4] = x_7 \quad [x_2, x_5] = x_7\quad [x_3, x_4] = -x_7$.\\
SQ.I. \quad $i = 1 \quad r=1 \quad c = 4 \quad Z = k[x_7]\\
F = <x_7>\quad CPI = <x_4, x_5, x_6, x_7>$\\
$T = <t>, \quad t = \diag (0,1,1,1,1,2,2)$.

\item[30.] $\mbox{\myfonts g}_{7,1.1(i_\lambda)} \cong (123457I) \cong R_{1}^\lambda, \lambda \neq 0,1$\\
$[x_1, x_2] = x_3 \quad [x_1, x_3] = x_4 \quad [x_1, x_4] = x_5\quad [x_1, x_5] = x_6 \quad [x_1, x_6] = x_7$\\
$[x_2, x_3] = x_5 \quad [x_2, x_4] = x_6 \quad [x_2, x_5] = \lambda x_7 \quad [x_3, x_4] = (1-\lambda) x_7$.\\
SQ.I. \quad $i = 1 \quad r=1 \quad c = 4 \quad Z = k[x_7]\\
F = <x_7>\quad CPI = <x_4, x_5, x_6, x_7>$\\
$T = <t>, \quad t = \diag (1,2,3,4,5,6,7)$.

\item[31.] ${\mbox{\myfonts g}}_{7,1.1(ii)} \cong (123457C) \cong R_{3}$\\
$[x_1, x_2] = x_3 \quad [x_1, x_3] = x_4 \quad [x_1, x_4] = x_5\quad [x_1, x_5] = x_6 \quad [x_1, x_6] = x_7$\\
$[x_2, x_5] = x_7 \quad [x_3, x_4] = -x_7$.\\
SQ.I. \quad $i = 1 \quad r=1 \quad c = 4 \quad Z = k[x_7]\\
F = <x_7>\quad CPI = <x_4, x_5, x_6, x_7>$\\
$T = <t>, \quad t = \diag (1,2,3,4,5,6,7)$.

\item[32.] ${\mbox{\myfonts g}}_{7,1.1(iv)} \cong (12457I) \cong R_{22}$\\
$[x_1, x_2] = x_3 \quad [x_1, x_3] = x_4 \quad [x_1, x_5] = x_6\quad [x_1, x_6] = x_7 \quad [x_2, x_3] = x_5$\\ 
$[x_2, x_4] = x_6 \quad [x_2, x_5] = x_7 \quad [x_3, x_4] = x_7$.\\
SQ.I. \quad $i = 1 \quad r=1 \quad c = 4 \quad Z = k[x_7]\\
F = <x_7>\quad CPI = <x_4, x_5, x_6, x_7>$\\
$T = <t>, \quad t = \diag (1,2,3,4,5,6,7)$.

\item[33.] ${\mbox{\myfonts g}}_{7,1.1(v)} \cong (12357C) \cong R_{42}$\\
$[x_1, x_3] = x_4 \quad [x_1, x_4] = x_5 \quad [x_1, x_5] = x_6\quad [x_1, x_6] = x_7 \quad [x_2, x_3] = x_5$\\
$[x_2, x_4] = x_6 \quad [x_2, x_5] = x_7 \quad [x_3, x_4] = -x_7$.\\
SQ.I. \quad $i = 1 \quad r=1 \quad c = 4 \quad Z = k[x_7]\\
F = <x_7>\quad CPI = <x_4, x_5, x_6, x_7>$\\
$T = <t>, \quad t = \diag (1,2,3,4,5,6,7)$.

\item[34.] ${\mbox{\myfonts g}}_{7,1.1(vi)} \cong (13457E) \cong R_{38}$\\
$[x_1, x_2] = x_3 \quad [x_1, x_3] = x_4 \quad [x_1, x_4] = x_5\quad [x_1, x_6] = x_7 \quad [x_2, x_3] = x_5$\\
$[x_2, x_5] = x_7 \quad [x_3, x_4] = -x_7$.\\
SQ.I. \quad $i = 1 \quad r=1 \quad c = 4 \quad Z = k[x_7]\\
F = <x_7>\quad CPI = <x_4, x_5, x_6, x_7>$\\
$T = <t>, \quad t = \diag (1,2,3,4,5,6,7)$.

\item[35.] ${\mbox{\myfonts g}}_{7,1.2(i_\lambda)} \cong (1357S)(\xi \neq 0) \cong R_{52}^\lambda, \lambda \neq 1$\\
$[x_1, x_2] = x_4 \quad [x_1, x_3] = x_6 \quad [x_1, x_4] = x_5\quad [x_1, x_5] = x_7 \quad [x_2, x_3] = \lambda x_5$\\
$[x_2, x_4] = x_6 \quad [x_2, x_6] = x_7 \quad [x_3, x_4] = (1-\lambda)x_7$.\\
SQ.I. \quad $i = 1 \quad r=1 \quad c = 4 \quad Z = k[x_7]\\
F = <x_7>\quad CPI = <x_4, x_5, x_6, x_7>$\\
$T = <t>, \quad t = \diag (1,1,2,2,3,3,4)$.

\item[36.] ${\mbox{\myfonts g}}_{7,1.2(ii)} \cong (1357S) (\xi = 0)$\\
$[x_1, x_2] = x_4 \quad [x_1, x_4] = x_5 \quad [x_1, x_5] = x_7\quad [x_1, x_6] = x_7 \quad [x_2, x_3] = x_6$\\
$[x_2, x_4] = x_6 \quad [x_2, x_5] = x_7 \quad [x_3, x_4] = -x_7$.\\
SQ.I. \quad $i = 1 \quad r=1 \quad c = 4 \quad Z = k[x_7]\\
F = <x_7>\quad CPI = <x_4, x_5, x_6, x_7>$\\
$T = <t>, \quad t = \diag (1,1,2,2,3,3,4)$.

\item[37.] ${\mbox{\myfonts g}}_{7,1.2(iv)} \cong (1357H) \cong R_{72}$\\
$[x_1, x_2] = x_4 \quad [x_1, x_4] = x_6 \quad [x_1, x_5] = -x_7\quad [x_1, x_6] = x_7 \quad [x_2, x_3] = x_5$\\
$[x_2, x_5] = x_7 \quad [x_3, x_4] = x_7$.\\
SQ.I. \quad $i = 1 \quad r=1 \quad c = 4 \quad Z = k[x_7]\\
F = <x_7>\quad CPI = <x_4, x_5, x_6, x_7>$\\
$T = <t>, \quad t = \diag (1,1,2,2,3,3,4)$.

\item[38.] ${\mbox{\myfonts g}}_{7,1.3(i_\lambda)} \cong (1357N) (\xi \neq 0) \cong R_{62}, \lambda \neq 0$\\
$[x_1, x_2] = x_4 \quad [x_1, x_3] = x_5 \quad [x_1, x_4] = x_6\quad [x_1, x_6] = x_7 \quad [x_2, x_3] = x_6$\\
$[x_2, x_4] = \lambda x_7 \quad [x_2, x_5] = x_7\quad [x_3,x_5] = x_7$.\\
SQ.I. \quad $i = 1 \quad r=1 \quad c = 4 \quad Z = k[x_7]\\
F = <x_7>\quad CPI = <x_4, x_5, x_6, x_7>$\\
$T = <t>,\quad t = \diag (1,2,2,3,3,4,5)$.

\item[39.] ${\mbox{\myfonts g}}_{7,1.3(ii)} \cong (1357L) \cong R_{63}$\\
$[x_1, x_2] = x_4 \quad [x_1, x_3] = x_5 \quad [x_1, x_4] = x_6\quad [x_1, x_6] = x_7 \quad [x_2, x_3] = x_6$\\
$[x_2, x_4] = x_7 \quad [x_2, x_5] = x_7/2 \quad [x_3,x_4] = -x_7/2$.\\
SQ.I. \quad $i = 1 \quad r=1 \quad c = 4 \quad Z = k[x_7]\\
F = <x_7>\quad CPI = <x_4, x_5, x_6, x_7>$\\
$T = <t>, \quad t = \diag (1,2,2,3,3,4,5)$.

\item[40.] ${\mbox{\myfonts g}}_{7,1.3(iii)} \cong (1357F) \cong R_{65}$\\
$[x_1, x_2] = x_4 \quad [x_1, x_3] = x_5 \quad [x_1, x_4] = x_6\quad [x_1, x_6] = x_7 \quad [x_2, x_4] = x_7$\\
$[x_3, x_5] = x_7$.\\
SQ.I. \quad $i = 1 \quad r=1 \quad c = 4 \quad Z = k[x_7]\\
F = <x_7>\quad CPI = <x_4, x_5, x_6, x_7>$\\
$T = <t>,\quad t = \diag (1,2,2,3,3,4,5)$.

\item[41.] ${\mbox{\myfonts g}}_{7,1.3(v)} \cong (1357C) \cong R_{99}$\\
$[x_1, x_2] = x_4 \quad [x_1, x_4] = x_6 \quad [x_1, x_6] = x_7\quad [x_2, x_3] = x_6 \quad [x_2, x_4] = x_7$\\
$[x_3, x_4] = -x_7 \quad [x_3, x_5] = -x_7$.\\
SQ.I. \quad $i = 1 \quad r=1 \quad c = 4 \quad Z = k[x_7]\\
F = <x_7>\quad CPI = <x_4, x_5, x_6, x_7>$\\
$T = <t>,\ t = \diag (1,2,2,3,3,4,5)$.

\item[42.] ${\mbox{\myfonts g}}_{7,1.8} \cong (1357J) \cong R_{70}$\\
$[x_1, x_2] = x_4 \quad [x_1, x_4] = x_6 \quad [x_1, x_6] = x_7\quad [x_2, x_3] = x_5 \quad [x_2, x_4] = x_7$\\
$[x_3, x_5] = x_7$.\\
SQ.I. \quad $i = 1 \quad r=1 \quad c = 4 \quad Z = k[x_7]\\
F = <x_7>\quad CPI = <x_4, x_5, x_6, x_7>$\\
$T = <t>,\quad t = \diag (2,4,3,6,7,8,10)$.\\

\item[43.] ${\mbox{\myfonts g}}_{7,1.11} \cong (12457E) \cong R_{46}$\\
$[x_1, x_2] = x_4 \quad [x_1, x_4] = x_5 \quad [x_1, x_5] = x_6\quad [x_1, x_6] = x_7 \quad [x_2, x_3] = x_6$\\ 
$[x_2, x_4] = x_6 \quad [x_2, x_5] = x_7 \quad [x_3,x_4] = -x_7$.\\
SQ.I. \quad $i = 1 \quad r=1 \quad c = 4 \quad Z = k[x_7]\\
F = <x_7>\quad CPI = <x_4, x_5, x_6, x_7>$\\
$T = <t>, \quad t = \diag (1,2,3,3,4,5,6)$.

\item[44.] ${\mbox{\myfonts g}}_{7,1.20} \cong (12457D) \cong R_{35}$\\
$[x_1, x_2] = x_3 \quad [x_1, x_5] = x_6 \quad [x_1, x_6] = x_7\quad [x_2, x_3] = x_4 \quad [x_2, x_4] = x_6$\\ 
$[x_2, x_5] = x_7 \quad [x_3, x_4] = x_7$.\\
SQ.I. \quad $i = 1 \quad r=1 \quad c = 4 \quad Z = k[x_7]\\
F = <x_7>\quad CPI = <x_4, x_5, x_6, x_7>$\\
$T = <t>, \quad t = \diag (1,2,3,5,6,7,8)$.

\item[45.] ${\mbox{\myfonts g}}_{7,2.1(i_\lambda)} \cong (1357M) \cong R_{64}^\lambda, \lambda \neq 0,1$\\
$[x_1, x_2] = x_4 \quad [x_1, x_3] = x_5 \quad [x_1, x_4] = x_6\quad [x_1, x_6] = x_7 \quad [x_2, x_3] = x_6$\\ 
$[x_2, x_5] = \lambda x_7 \quad [x_3, x_4] = (\lambda-1) x_7$.\\
SQ.I. \quad $i = 1 \quad r=2 \quad c = 4 \quad Z = k[x_7]\\
F = <x_7>\quad CPI = <x_4, x_5, x_6, x_7>$\\
$T = <t>, \quad t = \diag (\alpha, 1-3\alpha, 2\alpha, 1-2\alpha, 3\alpha, 1-\alpha, 1)$.

\item[46.] ${\mbox{\myfonts g}}_{7,2.1(ii)} \cong (1357D) \cong R_{69}$\\
$[x_1, x_2] = x_4 \quad [x_1, x_3] = x_5 \quad [x_1, x_4] = x_6\quad [x_1, x_6] = x_7 \quad [x_2, x_5] = x_7$\\ 
$[x_3, x_4] = x_7$.\\
SQ.I. \quad $i = 1 \quad r=2 \quad c = 4 \quad Z = k[x_7]\\
F = <x_7>\quad CPI = <x_4, x_5, x_6, x_7>$\\
$T = <t>, \quad t = \diag (\alpha, 1-3\alpha, 2\alpha, 1-2\alpha, 3\alpha,1-\alpha,1)$.

\item[47.] ${\mbox{\myfonts g}}_{7,2.1(iii)} \cong (1357A) \cong R_{101}$\\
$[x_1, x_2] = x_4 \quad [x_1, x_4] = x_6 \quad [x_1, x_6] = x_7\quad [x_2, x_3] = x_6 \quad [x_2, x_5] = x_7$\\ 
$[x_3, x_4] = -x_7$.\\
SQ.I. \quad $i = 1 \quad r=2 \quad c = 4 \quad Z = k[x_7]\\
F = <x_7>\quad CPI = <x_4, x_5, x_6, x_7>$\\
$T = <t>, \quad t = \diag (\alpha, 1-3\alpha, 2\alpha, 1-2\alpha, 3\alpha, 1-\alpha, 1)$.

\item[48.] ${\mbox{\myfonts g}}_{7,2.1(iv)} \cong (137D)$\\
$[x_1, x_3] = x_5 \quad [x_1, x_4] = x_6 \quad [x_1, x_6] = x_7\quad [x_2, x_3] = x_6 \quad [x_2, x_5] = x_7$\\ 
$[x_3, x_4] = x_7 $.\\
SQ.I. \quad $i = 1 \quad r=2 \quad c = 4 \quad Z = k[x_7]\\
F = <x_7>\quad CPI = <x_4, x_5, x_6, x_7>$\\
$T = <t>, \quad t = \diag (\alpha, 1-3\alpha, 2\alpha, 1-2\alpha, 3\alpha, 1-\alpha, 1)$.

\item[49.] ${\mbox{\myfonts g}}_{7,2.2} \cong (147D) \cong R_{95}$\\
$[x_1, x_2] = x_5 \quad [x_1, x_3] = x_6 \quad [x_1, x_4] = 2x_7\quad [x_2, x_3] = x_4 \quad [x_2, x_6] = x_7$\\ 
$[x_3, x_5] = -x_7 \quad [x_3,x_6] = x_7$.\\
SQ.I. \quad $i = 1 \quad r=2 \quad c = 4 \quad Z = k[x_7]\\
F = <x_7>\quad CPI = <x_4, x_5, x_6, x_7>$\\
$T = <t>, \quad t = \diag (1-2\alpha, \alpha, \alpha, 2\alpha, 1-\alpha, 1-\alpha, 1)$.

\item[50.] ${\mbox{\myfonts g}}_{7,2.10} \cong (13457C) \cong R_{39}$\\
$[x_1, x_2] = x_4 \quad [x_1, x_3] = x_7 \quad [x_1, x_4] = x_5\quad [x_1, x_5] = x_6 \quad [x_2, x_6] = x_7$\\ 
$[x_4, x_5] = -x_7 $.\\
SQ.I. \quad $i = 1 \quad r=2 \quad c = 4 \quad Z = k[x_7]\\
F = <x_7>\quad CPI = <x_3, x_5, x_6, x_7>$\\
$T = <t>, \quad t = \diag (-1+2\alpha, 2-3\alpha, 2-2\alpha, 1-\alpha, \alpha, -1+3\alpha, 1)$.

\item[51.] ${\mbox{\myfonts g}}_{7,2.23} \cong (137B)\cong R_{107}$\\
$[x_1, x_4] = x_6 \quad [x_1, x_6] = x_7 \quad [x_2, x_3] = x_5\quad [x_2, x_5] = x_7 \quad [x_3, x_4] = x_7$.\\
SQ.I. \quad $i = 1 \quad r=2 \quad c = 4 \quad Z = k[x_7]\\
F = <x_7>\quad CPI = <x_4, x_5, x_6, x_7>$\\
$T = <t>, \quad t = \diag (\alpha, \frac{1}{2} - \alpha, 2\alpha, 1-2\alpha, \frac{1}{2}+\alpha, 1-\alpha, 1)$.

\item[52.] ${\mbox{\myfonts g}}_{7,2.28} \cong (147B)\cong R_{112}$\\
$[x_1, x_2] = x_5 \quad [x_1, x_3] = x_6 \quad [x_1, x_6] = x_7\quad [x_2, x_5] = x_7 \quad [x_3, x_4] = x_7$.\\
SQ.I. \quad $i = 1 \quad r=2 \quad c = 4 \quad Z = k[x_7]\\
F = <x_7>\quad CPI = <x_4, x_5, x_6, x_7>$\\
$T = <t>, \quad t = \diag (1-2\alpha, \alpha, -1+4\alpha, 2-4\alpha, 1-\alpha, 2\alpha, 1)$.

\item[53.] ${\mbox{\myfonts g}}_{7,2.30} \cong (1457B)\cong R_{102}$\\
$[x_1, x_2] = x_5 \quad [x_1, x_5] = x_6 \quad [x_1, x_6] = x_7\quad [x_2, x_5] = x_7 \quad [x_3, x_4] = x_7$.\\
SQ.I. \quad $i = 1 \quad r=2 \quad c = 4 \quad Z = k[x_7]\\
F = <x_7>\quad CPI = <x_4, x_5, x_6, x_7>$\\
$T = <t>, \quad t = \diag (\frac{1}{5},\frac{2}{5}, \alpha, 1-\alpha, \frac{3}{5}, \frac{4}{5}, 1)$.

\item[54.] ${\mbox{\myfonts g}}_{7,2.37} \cong (1357R)$\\
$[x_1, x_2] = x_4 \quad [x_1, x_3] = x_5 \quad [x_1, x_4] = x_5\quad [x_1, x_6] = x_7 \quad [x_2, x_4] = x_6$\\
$[x_2, x_5] = x_7 \quad [x_3,x_4] = x_7$.\\
SQ.I. \quad $i = 1 \quad r=2 \quad c = 4 \quad Z = k[x_7]\\
F = <x_7>\quad CPI = <x_4, x_5, x_6, x_7>$\\
$T = <t>, \quad t = \diag (\alpha, \frac{1}{2}-\alpha, \frac{1}{2},\frac{1}{2},\frac{1}{2} +\alpha, 1-\alpha, 1)$.

\item[55.] ${\mbox{\myfonts g}}_{7,3.1(i_\lambda)} \cong (147E) \cong R_{93}^\lambda, \lambda \neq 0, 1$\\
$[x_1, x_2] = x_4 \quad [x_1, x_3] = x_5 \quad [x_1, x_6] = x_7\quad [x_2, x_3] = x_6 \quad [x_2, x_5] = \lambda x_7$\\ 
$[x_3, x_4] = (\lambda-1) x_7 $.\\
SQ.I. \quad $i = 1 \quad r=3 \quad c = 4 \quad Z = k[x_7]\\
F = <x_7>\quad CPI = <x_4, x_5, x_6, x_7>$\\
$T = <t>, \quad t = \diag (\alpha, \beta, 1-\alpha-\beta, \alpha+ \beta , 1-\beta, 1-\alpha, 1)$.

\item[56.] ${\mbox{\myfonts g}}_{7,3.1(i_\lambda)}, \lambda = 0$\\
$[x_1, x_2] = x_4 \quad [x_1, x_3] = x_5 \quad [x_1, x_6] = x_7\quad [x_2, x_3] = x_6 \quad [x_3, x_4] = -x_7$.\\
$i = 3 \quad r=3 \quad c = 5 \quad Z = k[x_5,x_7,x_2x_7-x_4x_6]\\
F = <x_2,x_4,x_5,x_6,x_7> = CPI $

\item[57.] ${\mbox{\myfonts g}}_{7,3.1(i_\lambda)}, \lambda = 1$\\
$[x_1, x_2] = x_4 \quad [x_1, x_3] = x_5 \quad [x_1, x_6] = x_7\quad [x_2, x_3] = x_6 \quad [x_2, x_5] = x_7$.\\
$i = 3 \quad r=3 \quad c = 5 \quad Z = k[x_4,x_7,x_3x_7-x_5x_6]\\
F = <x_3,x_4,x_5,x_6,x_7> = CPI $

\item[58.] ${\mbox{\myfonts g}}_{7,3.1(iii)} \cong (147A)\cong R_{113}$\\
$[x_1, x_2] = x_4 \quad [x_1, x_3] = x_5 \quad [x_1, x_6] = x_7\quad [x_2, x_5] = x_7 \quad [x_3, x_4] = x_7$.\\
SQ.I. \quad $i = 1 \quad r=3 \quad c = 4 \quad Z = k[x_7]\\
F = <x_7>\quad CPI = <x_4, x_5, x_6, x_7>$\\
$T = <t>, \quad t = \diag (\alpha, \beta, 1-\alpha-\beta,\alpha+\beta,1-\beta, 1-\alpha, 1)$.

\item[59.] ${\mbox{\myfonts g}}_{7,3.3} \cong (2457B) \cong R_{80}$\\
$[x_1, x_2] = x_4 \quad [x_1, x_4] = x_6 \quad [x_1, x_6] = x_7\quad [x_2, x_3] = x_5$.\\
$i = 3 \quad r=3 \quad c = 5 \quad Z = k[x_5, x_7, x_6^2-2x_4x_7]\\
F = <x_4, x_5, x_6, x_7>\quad CPI = <x_3, x_4, x_5, x_6, x_7>$\\
$T = <t_1, t_2, t_3>, \quad t_1 = \diag (0,1,0,1,1,1,1)$,\\
$t_2 = \diag (-1,3,0,2,3,1,0), \quad t_3 = \diag (0,0,1,0,1,0,0)$.

\item[60.] ${\mbox{\myfonts g}}_{7,3.4} \cong (247F) \cong R_{83}$\\
$[x_1, x_2] = x_4 \quad [x_1, x_3] = x_5 \quad [x_2, x_4] = x_6\quad [x_3, x_5] = x_7$.\\
$i = 3 \quad r=3 \quad c = 5 \quad Z = k[x_6, x_7, x_4^2x_7 +2x_1x_6x_7 + x_5^2x_6]\\
F = <x_1, x_4, x_5, x_6, x_7> = CPI$\\
$T = <t_1, t_2, t_3>, \quad t_1 = \diag (1,0,0,1,1,1,1)$,\\
$t_2 =\diag (0,1,0,1,0,2,0), \quad t_3 = \diag (2,-1,-1,1,1,0,0)$.

\item[61.] ${\mbox{\myfonts g}}_{7,3.5} \cong (247I) \cong R_{86}$\\
$[x_1, x_2] = x_4 \quad [x_1, x_3] = x_5 \quad [x_2, x_4] = x_6\quad [x_2, x_5] = x_7\quad [x_3,x_4] = x_7$.\\
$i = 3 \quad r=3 \quad c = 5 \quad Z = k[x_6, x_7, 2x_1x_7^2 + 2x_4x_5x_7 - x_5^2x_6]\\
F = <x_1, x_4, x_5, x_6, x_7> = CPI$\\
$T = <t_1, t_2, t_3>, \quad t_1 = \diag (1,0,0,1,1,1,1),\\
t_2 = \diag (1,-1,0,0,1,-1,0), \quad t_3 = \diag (2,-1,-1,1,1,0,0)$.

\item[62.] ${\mbox{\myfonts g}}_{7,3.6} \cong (357A) \cong R_{98}$\\
$[x_1, x_2] = x_4 \quad [x_1, x_3] = x_5 \quad [x_1, x_5] = x_7\quad [x_2, x_3] = x_6$.\\
SQ.I. \quad $i = 3 \quad r=3 \quad c = 5 \quad Z = k[x_4, x_6, x_7]\\
F = <x_4, x_6, x_7>\quad CPI = <x_3, x_4, x_5, x_6, x_7>$\\
$T = <t_1, t_2, t_3>, \quad t_1 = \diag (0,0,1,0,1,1,1),\\
t_2 = \diag (0,1,0,1,0,1,0), \quad t_3 = \diag (1,2,-2,3,-1,0,0)$.

\item[63.] ${\mbox{\myfonts g}}_{7,3.7} \cong (257H) \cong R_{117}$\\
$[x_1, x_2] = x_5 \quad [x_1, x_5] = x_6 \quad [x_2, x_4] = x_6\quad [x_3, x_4] = -x_7$.\\
$i = 3 \quad r=3 \quad c = 5 \quad Z = k[x_6, x_7, 2x_3x_6^2+2x_2x_6x_7-x_5^2x_7]\\
F = <x_2,x_3, x_5, x_6, x_7> = CPI $\\
$T = <t_1, t_2, t_3>, \quad t_1 = diag (0,0,1,0,0,0,1),\\
t_2 = \diag (0,1,0,0,1,1,0), \quad t_3 = \diag (1,-2,-2,2,-1,0,0)$.

\item[64.] ${\mbox{\myfonts g}}_{7,3.8} \cong (257A) \cong R_{121}$\\
$[x_1, x_2] = x_5 \quad [x_1, x_3] = x_6 \quad [x_1, x_5] = x_7\quad [x_2, x_4] = x_7$.\\
$i = 3 \quad r=3 \quad c = 5 \quad Z = k[x_6, x_7, x_3x_7 - x_5x_6]\\
F = <x_3, x_5, x_6, x_7>\quad CPI = <x_3, x_4, x_5, x_6, x_7>$\\
$T = <t_1, t_2, t_3>, \quad t_1 = \diag (1,0,0,2,1,1,2)\\
t_2 = \diag (0,0,1,0,0,1,0), \quad t_3 = \diag (1,-2,-1,2,-1,0,0)$.

\item[65.] ${\mbox{\myfonts g}}_{7,3.9} \cong (257C) \cong R_{123}$\\
$[x_1, x_2] = x_5 \quad [x_1, x_5] = -x_7 \quad [x_2, x_3] = x_6\quad [x_2, x_4] = x_7$.\\
$i = 3 \quad r=3 \quad c = 5 \quad Z = k[x_6, x_7, x_3x_7-x_4x_6]\\
F = <x_3,x_4, x_6, x_7>\quad CPI = <x_3, x_4, x_5, x_6, x_7>$\\
$T = <t_1, t_2, t_3>, \quad t_1 = \diag (0,1,0,0,1,1,1), \\
t_2 = \diag(0,0,1,0,0,1,0), \quad t_3 = \diag (1,-2,2,2,-1,0,0)$.

\item[66.] ${\mbox{\myfonts g}}_{7,3.10} \cong (137C) \cong R_{111}$\\
$[x_1, x_2] = x_5 \quad [x_1, x_3] = x_6 \quad [x_2, x_4] = x_6\quad [x_2, x_6] = x_7\quad [x_3,x_5] = x_7$.\\
$i = 3 \quad r=3 \quad c = 5 \quad Z = k[x_7, x_6^2 -2x_4x_7, x_1x_7 + x_5x_6]\\
F = <x_1, x_4, x_5, x_6, x_7> = CPI$\\
$T = <t_1, t_2, t_3>, \quad t_1 = \diag (1,0,0,1,1,1,1), \\
t_2 = \diag (0,-1,1,2,-1,1,0), \quad t_3 = \diag (1,0,-1,0,1,0,0)$.

\item[67.] ${\mbox{\myfonts g}}_{7,3.11} \cong (257B) \cong R_{119}$\\
$[x_1, x_2] = x_5 \quad [x_1, x_3] = x_6 \quad [x_1, x_5] = x_7\quad [x_2, x_4] = x_6$.\\
$i = 3 \quad r=3 \quad c = 5 \quad Z = k[x_6, x_7, x_3x_7 - x_5x_6]\\
F = <x_3, x_5, x_6, x_7>\quad CPI = <x_3, x_4, x_5, x_6, x_7>$\\
$T = <t_1, t_2, t_3>, \quad t_1 = \diag (0,1,0,-1,1,0,1), \\
t_2 = \diag (0,0,1,1,0,1,0), \quad t_3 = \diag (-1,2,1,-2,1,0,0)$.

\item[68.] ${\mbox{\myfonts g}}_{7,3.12} \cong (37D) \cong R_{124}$\\
$[x_1, x_2] = x_5 \quad [x_1, x_3] = x_6 \quad [x_2, x_4] = x_6\quad [x_3, x_4] = x_7$.\\
SQ.I. \quad $i = 3 \quad r=3 \quad c = 5 \quad Z = k[x_5, x_6, x_7]\\
F = <x_5, x_6, x_7>\quad CPI = <x_1, x_4, x_5, x_6, x_7>$.

\item[69.] ${\mbox{\myfonts g}}_{7,3.13} \cong (257K) \cong R_{105}$\\
$[x_1, x_2] = x_5 \quad [x_1, x_5] = x_6 \quad [x_2, x_5] = x_7\quad [x_3, x_4] = x_7$.\\
$i = 3 \quad r=3 \quad c = 5 \quad Z = k[x_6, x_7, x_5^2 + 2x_1x_7 - 2x_2x_6]\\
F = <x_1, x_2, x_5, x_6, x_7>$\quad No CP's.

\item[70.] ${\mbox{\myfonts g}}_{7,3.14} \cong (257F) \cong R_{120}$\\
$[x_1, x_2] = x_5 \quad [x_1, x_3] = x_6 \quad [x_1, x_6] = x_7\quad [x_2, x_4] = x_7$.\\
$i = 3 \quad r=3 \quad c = 5 \quad Z = k[x_5, x_7, x_6^2 -2x_3x_7]\\
F = <x_3, x_5, x_6, x_7>\quad CPI = <x_3, x_4, x_5, x_6, x_7>$\\
$T = <t_1, t_2, t_3>, \quad t_1 = \diag (0,0,1,1,0,1,1),\\
t_2 = \diag(-1,0,2,0,-1,1,0), \quad t_3 = \diag (0,1,0,-1,1,0,0)$.

\item[71.] ${\mbox{\myfonts g}}_{7,3.15} \cong (257E) \cong R_{122}$\\
$[x_1, x_2] = x_5 \quad [x_1, x_3] = x_6 \quad [x_2, x_5] = x_7\quad [x_3, x_4] = x_7$.\\
$i = 3 \quad r=3 \quad c = 5 \quad Z = k[x_6,x_7,x_5^2 + 2x_1x_7 + 2x_4x_6]\\
F = <x_1, x_4, x_5, x_6, x_7> = CPI$\\
$T = <t_1, t_2, t_3>, \quad t_1 = \diag (1,0,0,1,1,1,1),\\
t_2 = \diag (0,0,1,-1,0,1,0), \quad t_3 = \diag (2,-1,-2,2,1,0,0)$.

\item[72.] ${\mbox{\myfonts g}}_{7,3.16} \cong (137A) \cong R_{108}$\\
$[x_1, x_2] = x_5 \quad [x_1, x_5] = x_7 \quad [x_3, x_4] = x_6\quad [x_3, x_6] = x_7$.\\
$i = 3 \quad r=3 \quad c = 5 \quad Z = k[x_7, x_5^2 - 2x_2x_7, x_6^2 -2x_4x_7]\\
F = <x_2, x_4, x_5, x_6, x_7> = CPI $\\
$T = <t_1, t_2, t_3>,\quad t_1 = \diag (0,1,0,1,1,1,1), \\
t_2 = \diag (0,0,-1,2,0,1,0), \quad t_3 = \diag (1,-2,0,0,-1,0,0)$.

\item[73.] ${\mbox{\myfonts g}}_{7,3.17} \cong (1457A) \cong R_{103}$\\
$[x_1, x_2] = x_5 \quad [x_1, x_5] = x_6 \quad [x_1, x_6] = x_7 \quad [x_3, x_4] = x_7$.\\
$i = 3 \quad r=3 \quad c = 5 \quad F = <x_2, x_5,x_6, x_7>$\\
$CPI = <x_2,x_4,x_5,x_6,x_7> \quad Z = k[x_7, f_1, f_2, f_3]$\\
$f_1 = 2x_5x_7 - x_6^2,\\
f_2 = 3x_2x_7^2 -3x_5x_6x_7 + x_6^3,\\
f_3 = 9x_2^2x_7^2 - 18x_2x_5x_6x_7 + 6x_2x_6^3 + 8x_5^3x_7 - 3x_5^2x_6^2$.\\
Relation: \quad $f_1^3 + f_2^2 - x_7^2f_3 = 0\\
Q(Z) = k(x_7, f_1, f_2)$. \quad Note that $Z$ is isomorphic with $Z(U({\mbox{\myfonts g}}_{5,5}))$.

\item[74.] ${\mbox{\myfonts g}}_{7,3.18} \cong (157) \cong R_{129}$\\
$[x_1, x_2] = x_6 \quad [x_1, x_6] = x_7 \quad [x_2, x_5] = x_7\quad [x_3, x_4] = x_7$.\\
SQ.I. \quad $i = 1 \quad r=3 \quad c = 4 \quad Z = k[x_7]\\
F = <x_7> \quad CPI = <x_4, x_5, x_6, x_7>$\\
$T = <t>, \quad t = \diag (\alpha, 1-2\alpha, \beta,1-\beta, 2\alpha, 1-\alpha, 1)$.

\item[75.] ${\mbox{\myfonts g}}_{7,3.19} \cong (27B) \cong R_{130}$\\
$[x_1, x_2] = x_6 \quad [x_1, x_3] = x_7 \quad [x_3, x_4] = x_6\quad [x_4, x_5] = x_7$.\\
$i = 3 \quad r=3 \quad c = 5 \quad Z = k[x_6, x_7, x_2x_7^2 - x_3x_6x_7 -x_5x_6^2]\\
F = <x_2, x_3, x_5, x_6, x_7> = CPI $\\
$T = <t_1, t_2, t_3>, \quad t_1 = \diag (0,0,1,-1,2,0,1), \\
t_2 = \diag (0,1,0,1,-1,1,0), \quad t_3 = \diag (1,-1,-1,1,-1,0,0)$.

\item[76.] ${\mbox{\myfonts g}}_{7,3.21} \cong (247B)$\\
$[x_1, x_2] = x_4 \quad [x_1, x_3] = x_5 \quad [x_1, x_4] = x_6\quad [x_3, x_5] = x_7$.\\
$i = 3 \quad r=3 \quad c = 5 \quad Z = k[x_6, x_7, x_4^2 - 2x_2x_6]\\
F = <x_2, x_4, x_6, x_7>\quad CPI = <x_2, x_4, x_5, x_6, x_7>$\\
$T = <t_1, t_2, t_3>, \quad t_1 = \diag (1,0,0,1,1,2,1), \\
t_2 = \diag (0,1,0,1,0,1,0), \quad t_3 = \diag (2,-4,-1,-2,1,0,0)$.

\item[77.] ${\mbox{\myfonts g}}_{7,3.22} \cong (247D)$\\
$[x_1, x_2] = x_4 \quad [x_1, x_3] = x_5 \quad [x_1, x_5] = x_7\quad [x_2, x_5] = x_6\quad [x_3, x_4] = x_6$.\\
$i = 3 \quad r=3 \quad c = 5 \quad Z = k[x_6, x_7, x_1x_6 - x_2x_7+x_4x_5]\\
F = <x_1, x_2, x_4, x_5, x_6, x_7>$\quad No CP's \\
$T = <t_1, t_2, t_3>, \quad t_1 = \diag (0,0,1,0,1,1,1,1), \\
t_2 = \diag (0,1,0,1,0,1,0), \quad t_3 = \diag(1,1,-2,2,-1,0,0)$.

\item[78.] ${\mbox{\myfonts g}}_{7,3.23} \cong (357B)$\\
$[x_1, x_2] = x_3 \quad [x_1, x_3] = x_5 \quad [x_1, x_4] = x_7\quad [x_2, x_3] = x_6$.\\
SQ.I. $i = 3 \quad r=3 \quad c = 5 \quad Z = k[x_5, x_6, x_7]\\
F = <x_5, x_6, x_7>\quad CPI = <x_3, x_4, x_5, x_6, x_7>$\\
$T = <t_1, t_2, t_3>, \quad t_1 = \diag (0,0,0,1,0,0,1), \\
t_2 = \diag(0,1,1,0,1,2,0), \quad t_3 = \diag(2,-1,1,-2,3,0,0)$.\\

\item[79.] ${\mbox{\myfonts g}}_{7,3.24} \cong (37C)$\\
$[x_1, x_2] = x_5 \quad [x_2, x_3] = x_6 \quad [x_2, x_4] = x_7\quad [x_3, x_4] = x_5$.\\
SQ.I. \quad $i = 3 \quad r=3 \quad c = 5 \quad Z = k[x_5, x_6, x_7]\\
F = <x_5, x_6, x_7>\quad CPI = <x_1, x_4, x_5, x_6, x_7>$\\
$T = <t_1, t_2, t_3>, \quad t_1 = \diag (0,0,-1,1,0,-1,1),\\
t_2 = \diag (1,0,1,0,1,1,0), \quad t_3 = \diag(3,-1,1,1,2,0,0)$.

\item[80.] ${\mbox{\myfonts g}}_{7,4.1} \cong (37B) \cong R_{126}$\\
$[x_1, x_2] = x_5 \quad [x_1, x_3] = x_6 \quad [x_3, x_4] = x_7$.\\
SQ.I. \quad $i = 3 \quad r=4 \quad c = 5 \quad Z = k[x_5, x_6, x_7]\\
F = <x_5, x_6, x_7>\quad CPI = <x_2, x_4, x_5, x_6, x_7>$\\
$T = <t_1, t_2, t_3>, \quad t_1 = \diag (\alpha, 1-\alpha, 1-\alpha, \alpha, 1,1,1), \\
t_2 = \diag (1,0,0,0,1,1,0), \quad t_3 = \diag (1,0,-1,1,1,0,0)$.

\item[81.] ${\mbox{\myfonts g}}_{7,4.3} \cong (27A) \cong R_{131}$\\
$[x_1, x_2] = x_6 \quad [x_3, x_5] = x_6 \quad [x_4, x_5] = x_7$.\\
$i = 3 \quad r=4 \quad c = 5 \quad Z = k[x_6, x_7, x_3x_7 - x_4x_6]\\
F = <x_3, x_4, x_6, x_7>\quad CPI = <x_1, x_3, x_4, x_6, x_7>$\\
$T = <t_1, t_2, t_3>, \quad t_1 = \diag (\alpha, 1-\alpha, 1, 1, 0,1,1),\\
t_2 = \diag (1,0,1,0,0,1,0), \quad t_3 = \diag (0,0,1,1,-1,0,0)$.

\item[82.] ${\mbox{\myfonts g}}_{7,4.4} \cong (17) \cong R_{132}$ (7-dim Heisenberg Lie algebra)\\
$[x_1, x_4] = x_7 \quad [x_2, x_5] = x_7 \quad [x_3, x_6] = x_7$.\\
SQ.I. \quad $i = 1 \quad r=4 \quad c = 4 \quad Z = k[x_7]\\
F = <x_7> \quad CPI = <x_4, x_5, x_6, x_7>$\\
$T = <t>, \quad t = \diag (\alpha, \beta, \gamma, 1-\alpha, 1-\beta, 1-\gamma, 1)$.
\end{itemize}
\ \\
\begin{center}
\ \\
{\bf ACKNOWLEDGMENTS}
\end{center}
We like to thank Alexander Elashvili for fruitful discussions on the subject and for his encouragement.  We are also very grateful to Michel Van den Bergh for his genuine interest and for his valuable help in calculating the invariants for $\mbox{\myfonts g}_{6,6}$ and $\mbox{\myfonts g}_{6,16}$.  Finally, we wish to thank Mustapha Rais and Fran\c{c}ois Rouvi\`ere for sending us some unpublished manuscripts by Andr\'e Cerezo (1945-2003).  The life and work of this remarkable mathematician is commemorated in the interesting website 
\begin{center}
http://math.unice.fr/$\sim$frou/AC.html
\end{center}
created by Fran\c{c}ois Rouvi\`ere.

\end{document}